\documentclass{amsart}
\usepackage{amssymb}  
\usepackage[all]{xy}  
  
\newcommand\bS{{}^b\kern-2ptS}  
\newcommand\PsS[2]{\Psi_{\rm inv}^{#1}(#2)}  
\newcommand\Psm[2]{\psi_{\rm inv}^{#1}(#2)}  
\newcommand\Ps[2]{\Psi^{#1}(#2)}  
  
\newcommand\Hd{\operatorname{HH}}  
  
\newcommand\Hom{\operatorname{Hom}}  
\newcommand\alg{\operatorname{alg}}  
\newcommand\topo{\operatorname{top}}  
  
\newcommand\ord{\operatorname{ord}}  
  
\newcommand\aTr{\operatorname{\overline{Tr}}}  
\newcommand\ITr{\,\operatorname{pv}\!\!-\!\!\!\int}  
  
\newcommand\Tr{\operatorname{Tr}}  
\newcommand\bTr[1]{\overline{\operatorname{Tr}}_{#1}}  
\newcommand\RTr{\operatorname{Tr_R}}  
\newcommand\RT{\operatorname{Tr_R}}  
  
\newcommand\GR{\mathcal{G}}  
  
\newcommand\pSS[1]{{}^p\kern-1pt\SS{#1}}  
  
\newcommand\ind{\operatorname{ind}}   
\newcommand\mfk{\mathfrak}  
  
\newcommand\datver[1]{\def\datverp%
 {\par\boxed{\boxed{\text{Version: #1; Run: \today}}}}}  
\datver{3.2; Revised: July 31, 2000}  
  
\renewcommand\Re{\operatorname{Re}}

\newcommand\CC{\mathbb C}  
\newcommand\NN{\mathbb N}  
\newcommand\QQ{\mathbb Q}  
\newcommand\RR{\mathbb R}  
\newcommand\ZZ{\mathbb Z}  
\renewcommand\SS{\mathbb S}  
  
\newcommand\lgg{\mathfrak g}  
\newcommand\pa{\partial}

\newcommand\CI{{\mathcal C}^{\infty}}  
  
\newcommand\CIc{{\mathcal C}^{\infty}_{\text{c}}}  
  
\newcommand\Id{\operatorname{Id}}  
\newcommand\Mand{\text{ and }}  
  
\newcommand\Mfor{\text{ for }}  
\newcommand\Mif{\text{ if }}  
  
\newcommand\ie{{\em i.e.,} }  
\newcommand\hotimes{\widehat\otimes}

\newtheorem{theorem}{Theorem}  
\newtheorem{proposition}{Proposition}  
\newtheorem{corollary}{Corollary}  
\newtheorem{lemma}{Lemma}  
  
\theoremstyle{definition}

\theoremstyle{remark}

\begin{document}  
\title[Index for families] {An index theorem for families invariant  
with respect to a bundle of Lie groups}  
  
\author{Victor Nistor} \address{Current address:\  
  Department of Mathematics,  
  Pennsylvania State University} \email{nistor@math.psu.edu}  
\thanks{Partially supported by the NSF Young Investigator Award  
  DMS-9457859, a Sloan fellowship, and NSF Grant  
DMS-9971951. {\bf  
    http:{\scriptsize//}www.math.psu.edu{\scriptsize/}nistor{\scriptsize/}}.  
}  
\dedicatory\datverp  

\begin{abstract} 
We define the equivariant family index of a family of elliptic
operators invariant with respect to the free action of a bundle $\GR$
of Lie groups. If the fibers of $\GR \to B$ are simply-connected
solvable, we then compute the Chern character of the (equivariant
family) index, the result being given by an Atiyah-Singer type
formula. We also study traces on the corresponding algebras of
pseudodifferential operators and obtain a local index formula for such
families of invariant operators, using the Fedosov product. For
topologically non-trivial bundles we have to use methods of
non-commutative geometry. We discuss then as an application the
construction of ``higher-eta invariants,'' which are morphisms
$K_n(\PsS {\infty}Y) \to \CC$.  The algebras of invariant
pseudodifferential operators that we study, $\Psm {\infty}Y$ and $\PsS
{\infty}Y$, are generalizations of ``parameter dependent'' algebras of
pseudodifferential operators (with parameter in $\RR^q$), so our
results provide also an index theorem for elliptic, parameter
dependent pseudodifferential operators.
\end{abstract}  
\maketitle  
  
\tableofcontents  
  
\section*{Introduction}

Families of Dirac operators invariant with respect to a bundle of  
Lie groups appear in the analysis of the Dirac operator on  
certain non-compact manifolds. They arise, for example, in the  
analysis of the Dirac operator on an $S^1$-manifold $M$, if we  
desingularize the action of $S^1$ by replacing the original metric $g$  
with $\phi^{-2}g$, where $\phi$ is the length of the infinitesimal  
generator $X$ of the $S^1$-action. In this way, $X$ becomes of length  
one in the new metric.  The main result of \cite{NistorDOP} states  
that the kernel of the new Dirac operator on the open manifold  
$M \smallsetminus M^{S^1}$ is naturally isomorphic to the kernel of the  
original Dirac operator.  
  
It turns out that the Fredholm property of the resulting Dirac  
operator (obtained by the above procedure on the non-compact manifold  
$M \smallsetminus M^{S^1}$) is controlled by the invertibility of a  
family of operators invariant with respect to the action of a bundle  
of non-abelian, solvable Lie groups. This follows from the results of  
\cite{LauterNistor} and it will be discussed in greater detail in a  
future paper. In general, neither the bundle $Y \to B$ on which these  
operators act, nor the bundle of Lie groups $\GR \to B$ acting on  
$Y$, are trivial. A natural problem then is to study the invertibility  
of these invariant families of operators, and more generally, their  
index.  
  
We define the (equivariant family) index of a family of invariant,  
elliptic operators using $K$-theory. For operators acting between  
sections of the same vector bundle, we can define the index using the  
boundary map in algebraic $K$-Theory. For operators acting between  
sections of different bundles, one has to use Kasparov's bivariant  
$K$-Theory \cite{Blackadar} or a smooth variant of it  
\cite{NistorKthry}. The algebraic $K$-theory definition of the index  
gives a little bit more than the one using bivariant $K$-theory, but  
it applies only to elliptic operators acting between {\em isomorphic}  
vector bundles, which is however almost always the case in  
applications, and hence only a minor drawback. In any case, it turns  
out that the (equivariant family) index of such an invariant elliptic  
family is the obstruction to finding an invertible perturbation of the  
original family by families of invariant, regularizing operators, if  
we exclude the degenerate case $\dim Y = \dim \GR$. This shows the  
relevance of computing the index to the problem of determining the  
invertibility of a given family.  
  
In this paper, we study the index and certain non-local invariants of  
families of elliptic operators invariant with respect to a bundle of  
{\em simply-connected, solvable} Lie groups $\GR$. One of the main  
results is a formula for the Chern character of the index bundle that  
is similar to the Atiyah-Singer index formula for families. We  
analyse the the local behavior of these families, when  
$\GR$ is a vector bundle. This leads us then to the construction  
of several traces on the algebras $\PsS {m}Y$. We use these traces to  
obtain local index theorems. For more general bundles $\GR$, the  
local analysis is likely to be much harder because there are no good   
candidates for the construction of convolution algebras   
closed under functional analytic calculus on non-commutative solvable   
Lie groups.   
  
In \cite{BismutCheeger}, Bismut and Cheeger have generalized the  
Atiyah-Patodi-Singer index theorem \cite{APS} to families of Dirac  
operators on manifolds with boundary (see also  
\cite{Melrose-Piazza1}). Their results apply to operators whose  
``indicial parts'' are invertible. These indicial parts are actually  
families of Dirac operators invariant with respect to a one-parameter  
group, so they fit into the framework of this paper (with $\GR = B  
\times \RR$).  In addition to the usual ingredients of an index  
theorem -- curvature and characteristic classes -- their result was 
stated in terms of a new invariant, called the ``eta-form,'' in 
analogy to the additional invariant appearing in the 
Atiyah-Patodi-Singer index formula for operators on manifolds with 
boundary. Thus, the results of this paper are relevant also to the 
problem of determining an explicit formula for the index of a family 
of pseudodifferential elliptic operators on a bundle of manifolds with 
boundary. (See also \cite{troitsky}.) 
  
With an eye towards this problem, we also give a new proof of the  
regularity of the eta function at the origin and discuss some possible  
generalizations the eta invariant. Actually, we discuss two possible  
generalizations, one which is a direct generalization of a result of  
\cite{Melrose46} and one using higher algebraic $K$-theory. The first  
possible generalization is to associate to a Dirac operator invariant  
with respect to $\RR^q$ the quantity defined by the formula $\bTr q  
\big( (D^{-1}dD)^{2k-1} \big)$.  This possible generalization was  
considered before by Lesch and Pflaum \cite{LP}, who proved that this  
formula does not lead to new invariants for Dirac operators and also  
that this formula is not additive for a product of two invertible  
operators, except for $k =1$, when one recovers the usual eta  
invariant \cite{Melrose46}. The second possible generalization, which  
has the advantage of being additive, is to define the higher eta  
invariant as a morphism on higher algebraic $K$-theory. This  
definition can be found in Section \ref{S.HEF}.  
  
We now describe the contents of each section of this paper. In Section  
\ref{S.Inv.Ops}, we discuss the action of a bundle of Lie groups $\GR$  
on a fiber bundle $Y$ and we introduce the algebras $\Psm {\infty}Y$  
and $\PsS {\infty}Y$, which will be our main object of study. (Both  
these algebras consist of families of invariant pseudodifferential  
operators.)  Families of operators invariant with respect to a bundle  
of Lie groups were probably considered for the first time in  
\cite{MelroseCongress}. (See also \cite{MazzeoMelrose}.) When $B$ is  
reduced to a point and $\GR = \RR^q$, the algebras $\PsS{\infty}Y$  
were studied in \cite{Melrose42}. Torsion invariants for families of  
Dirac operators invariant with respect to a vector bundle were defined  
and studied by Bismut, see \cite{BismutSURV}. We prove that, when  
$\GR$ is a vector bundle, the group of gauge transformations of $\GR$  
acts on $\Psm {\infty}Y$ and $\PsS {\infty}Y$. In Section  
\ref{S.Homotopy.Inv}, we define the index of a family of elliptic,  
invariant pseudodifferential operators $A$. We shall sometimes use the  
term ``the equivariant family index'' of an elliptic $A \in \Psm  
{\infty}Y$, for the index of such a family. We prove that the index of  
$A$ is the obstruction to finding a regularizing $R$ such that $A +  
R$, acting between suitable Sobolev spaces, is invertible in each  
fiber (excluding the degenerate case $\dim Y = \dim \GR$). This  
generalizes the usual property of the Fredholm index of a Fredholm  
operator.  In Section \ref{Sec.Solvable} we prove that when $\GR$  
consists of simply-connected solvable Lie groups  
\begin{equation}\label{eq.K.solv}  
	K_*(C^*(\GR)) \simeq K^*(\lgg),  
\end{equation}  
if $\lgg$ is the vector bundle of the Lie algebras defined  
by $\GR$. Then we obtain a formula for the Chern character of the  
index of an elliptic operator $A \in \Psm {\infty}Y$. (For simplicity,  
we called $A$ ``an operator,'' although it is really a family of  
operators. We shall do this repeatedly.)   
  
Beginning with Section \ref{S.RT}, we assume that $\GR$ is a vector  
bundle, because then we can construct natural algebras that are  
spectraly invariant (even closed under functional analytic calculus),  
and hence we obtain more refined results. We then develop the  
necessary facts about the asymptotics of the trace of the operators in  
$\PsS {\infty}Y$. This allows us to define various regularized and  
residue traces. Using these traces, we obtain in Section \ref{S.HI}  
two local formulae for the (equivariant family) index of an elliptic  
operator $A \in \PsS{\infty}Y$. In Section \ref{S.HEF},   
we discuss two  
possible generalizations of the eta invariant,   
which is suggested by the formula  
$$  
	\eta(D_0,s) = \bTr 1( D^{-1} D'), \quad \text{where } D = D_0 +  
	\pa_t \text{ and } D' = [D,t],  
$$  
proved in \cite{Melrose46} using the local  
index theorem (see also \cite{LP}). Operators invariant with respect  
to $\RR^q$, a particular case of our operators when the base is  
reduced to a point, appear in the formulation of elliptic (or  
Fredholm) boundary conditions for pseudodifferential operators on  
manifolds with corners. The equivariant index can be used to study  
this problem, which is relevant to the question of extending the  
Atiyah-Patodi-Singer boundary conditions to manifolds with corners.  
  
The algebras of invariant pseudodifferential operators that we study,  
$\Psm {\infty}Y$ and $\PsS {\infty}Y$, are generalizations of  
``parameter dependent'' algebras of pseudodifferential operators  
considered by Agmon \cite{Agmon}, Grubb and Seeley \cite{GrubbSeeley},  
Lesch and Pflaum \cite{LP}, Melrose \cite{Melrose46}, Shubin  
\cite{SH}, and others. Our index theorem, Theorem \ref{Theorem.Chern},  
hence gives a solution to the problem of determining the index of  
elliptic, parameter dependent families of pseudodifferential  
operators, parameterized by $\lambda \in \RR^q$. We note that the  
concept of index of such a family requires a proper definition, and  
that the Fredholm index, ``dimension of the kernel'' - ``dimension of  
the cokernel,'' is not appropriate for $q \ge 1$.  Our definition of  
the (equivariant, family) index is somewhat closer to the definition  
of the $L^2$-index for covering spaces given in \cite{Atiyah} and  
\cite{Singer} than to the definition of the Fredholm index. However,  
an essential difference between their definition and ours is that no  
trace is involved in our definition.  
  
For the proof of our  
local index theorem, we use ideas of non-commutative geometry  
\cite{ConnesNCG}, more precisely, the general approach to index  
theorems using cyclic cohomology as developed in \cite{Nistor2}. These  
computations are also an example of a computation of a bivariant  
Chern-Connes character \cite{NistorANN}.  We plan to use the results of  
this paper for some problems in $M$-theory \cite{SethiStern}. We also expect   
our results to have applications to adiabatic limits of eta invariants  
\cite{BF,Witten1}.  
  
There exist several potential extensions of our results, although none  
of them seems to be straightforward. These extensions will presumable  
involve more general conditions on $\GR$ and its action on  
$Y$. Probably the most general conditions under which one can  
reasonably expect to obtain definite results are those when the fibers  
of $\GR$ are connected and the action is proper. We assume that the  
fibers of $\GR$ are connected to be able to use general results on  
connected Lie groups.  
  
Therefore, let us consider a general family of connected Lie groups  
$\GR$ and successively weaken our assumptions on $\GR$ and its action  
on $Y$ and then try to see to what extend our results extend to this  
general setting.  If we continue to assume that the action of $\GR$ is  
free, then the results of Sections \ref{S.Inv.Ops} and  
\ref{S.Homotopy.Inv} remain true, because we do not use any conditions  
on the fibers of $\GR$ in these two section. However, the results of  
the subsequent sections are no longer valid.  If we drop the  
assumption that the action of $\GR$ be free, then the definition of  
the algebras $\Psm{\infty}Y$ still makes sense in this more general  
situation, provided that we assume the action of $\GR$ to be proper,  
or otherwise our algebras might be too small.  
  
Let us then assume henceforth that the action of $\GR$ is proper. For  
simply-connected solvable families, this implies that the action is  
free, which justifies our choices and moreover gives that $Y \cong Z  
\times_B \GR$, as $\GR$-spaces. For general families though, the  
topology on $Y$ and the action of $\GR$ will be more complicated if  
the action of $\GR$ is not free. Moreover, except the definition of  
the algebras $\Psm{\infty}Y$, few of our results extend to the general  
case of connected fibers and proper actions.  For example, the results  
of Section \ref{S.Homotopy.Inv} are not true for non-free actions, in  
general. This is the case, in particular, of Lemma \ref{lemma.stable},  
so we cannot define the equivariant family index as an element of  
$K_0(C^*_r(\GR))$, as we do in this paper, but we get it as an element  
of $K_0(C^*_r(Y;\GR))$. Further work is necessary to show that the  
natural map $K_0(C^*_r(Y;\GR)) \to K_0(C^*_r(\GR))$ is injective in  
the cases of interest.  
  
Let us now take a closer look at the case when the fibers of $\GR \to  
B$ are compact. Then we do get that $K_0(C^*_r(Y;\GR)) \to  
K_0(C^*_r(\GR))$ is injective, and hence we can still define the index  
as an element of $K_0(C^*_r(\GR))$.  It is not difficult to show that  
all the fibers $\GR_b$ are isomorphic as groups to a fixed Lie group,  
say $G$.  The fundamental group $\pi$ of $B$ then acts by holonomy on  
$R(G)$, the representation ring of $G$, and if $G$ is connected and  
semisimple, we have  
\begin{equation}  
	K_*(C^*(\GR)) \otimes \QQ \cong K^*(B) \otimes R(G)^{\pi}  
	\otimes \QQ,  
\end{equation}  
In general, however, there will be no such isomorphism if we do not  
tensor by $\QQ$.  Thus, when the fibers of $\GR$ are compact, the  
groups $K_*(C^*_r(\GR))$ are completely different from the  
corresponding groups when the fibers of $\GR$ are solvable (see also  
Equation \eqref{eq.K.solv}).  
  
Nevertheless, the methods of this paper will probably be useful to  
treat the general case, allowing us to reduce the general case of  
connected fibers and proper actions to the case of a bundle of  
semisimple Lie groups. The ``Dirac'' and ``dual Dirac'' bivariant  
$K$-theory elements will probably then allow us to further reduce the  
semisimple case to the compact case, in the spirit of the  
Connes-Kasparov conjecture, using Kasparov's bivariant $K$-theory. The  
case of compact fibers can then be treated directly, although this is  
not straigtforward and may require fields of matrices with non-trivial  
Diximier-Douady invariants.  Moreover, some conceptual difficulties in  
determining the equivariant index arise when the bundle $\GR$ is  
non-trivial.  
  
I would like to thank Fourier Institute, where part of this project  
was completed, for their hospitality. Also, I would like to thank  
Florian Albers, Jean-Michel Bismut, Alain Connes, David Kazhdan, and   
Richard Melrose for useful comments.  
  
All pseudodifferential operators considered in this paper are  
``classical,'' that is, one-step polyhomogeneous.

\section{Invariant pseudodifferential operators\label{S.Inv.Ops}}

We now describe the settings in which we shall work. Let $B$ be a  
smooth compact manifold and  
$$  
        d : \GR \to B \quad \mbox{ and } \pi : Y \to B  
$$  
be two smooth fiber bundles with fibers $\GR_b := d^{-1}(b)$ and $Y_b  
:= \pi^{-1}(b)$. We shall assume that $\GR$ is a bundle of Lie groups  
acting smoothly on $Y$, and then we shall consider families of  
operators along the fibers of $Y$ and invariant under the action of  
$\GR$. We can restrict in our discussion to a connected component of  
$B$, so for simplicity, we shall assume that $B$ is connected.  The  
index and local invariant of these operators will form our main object  
of study. We now make all these assumptions and concepts more precise.  
  
Throughout this paper, $\GR$ will denote a {\em bundle of Lie groups}  
on a manifold usually denoted $B$. By this we mean that $\GR \to B$ is  
a smooth fiber bundle, that each $\GR_b$ is a Lie group, and that the  
multiplication and inverse depend differentiably on $b$. Hence the map  
\begin{equation}  
        \GR \times_B \GR := \{ (g',g) \in \GR \times \GR, d(g') = d(g)\}  
        \ni\, (g',g) \longrightarrow g'g^{-1}\, \in \GR  
\end{equation}  
is differentiable. This implies, by standard arguments, that the map  
sending a point $b \in \GR_b$ to $e_b$, the identity element of  
$\GR_b$, is a diffeomorphism onto a smooth submanifold of $\GR$. It  
also implies that the map $\GR \in g \to g^{-1} \in \GR$ is  
differentiable.  
  
We also assume that $\GR$ acts smoothly on $Y$. This means that there  
are given actions $\GR_b \times Y_b \to Y_b$ of $\GR_b$ on $Y_b$, for  
each $b$, such that the induced map,  
$$  
        \GR \times_B Y :=\{ (g,y) \in \GR \times Y, d(g) = \pi(y)\}  
        \ni\, (g,y) \longrightarrow gy\, \in Y,  
$$  
is differentiable. We shall also assume that the action of $\GR$ on  
$Y$ is {\em free}, that is, that the action of $\GR_b$ on $Y_b$ is  
free for each $b$.  
  
We do not assume, however, that the groups $\GR_b$ are isomorphic,  
although this is true in most applications.  For us, two important  
particular cases of Lie groups bundles are when the fibers of $\GR$  
are compact and the case when the fibers of $\GR$ are  
simply-connected, solvable Lie groups. Let $\lgg \to B$ be the bundle  
of Lie algebras of associated to $\GR$, that is, the bundle whose  
fiber above $b$, $\lgg_b$, is $Lie(\GR_b)$. We shall also treat in  
detail the case when $\GR$ consists of simply-connected abelian  
groups, that is, it is a vector bundle $\GR \simeq \lgg$. This case is  
very important in applications and can be treated more completely  
because it does not involve any complications related to the harmonic  
analysis on non-compact, non-abelian Lie groups.  
  
On $Y$, we consider smooth families $A = (A_b)$, $b \in B$, of  
classical pseudodifferential operators acting on the fibers of $Y \to  
B$ such that each $A_b$ is invariant with respect to the action of the  
group $\GR_b$. Unless mentioned otherwise, we assume that these  
operators act on half densities along each fiber.  The algebra that we  
are interested in consists of such invariant operators satisfying also  
a {\em support condition}. To state this support condition, first  
notice that a family $A = (A_b)$ defines a continuous map $\CIc(Y) \to  
\CI(Y)$, and, as such, it has a distribution (or Schwartz) kernel,  
which is a distribution $K_A$ on $Y \times_B Y \subset Y \times Y$.  
(We ignore the vector bundles in which this distribution takes its  
values.)  Because the family $A=(A_b)$ is invariant, the distribution  
$K_A$ is also invariant with respect to the action of  
$\GR$. Consequently, $K_A$ is the pull back of a distribution $k_A$ on  
$(Y \times_B Y)/\GR$. {\em We will require that $k_A$ have compact  
support.} We shall sometimes call $k_A$ the {\em convolution kernel}  
of $A$.  This condition on the support of $k_A$ ensures that each  
$A_b$ is a properly supported pseudodifferential operator, and hence  
it maps compactly supported functions (or sections of a vector bundle,  
if we consider operators acting on sections of a smooth vector bundle)  
to compactly supported functions (or sections). This support condition  
is automatically satisfied if $Y/\GR$ is compact and each $A_b$ is a  
differential operator. The space of smooth, invariant families $A$ of  
order $m$ pseudodifferential operators acting on the fibers of $Y \to  
B$ such that $k_A$ has compact support will be denoted by $\Psm{m}Y$.  
Then  
$$  
	\Psm {\infty}Y := \cup_{m \in \ZZ} \Psm{m}Y  
$$   
is an algebra, by classical results \cite{Hormander3}. Note also that $\Psm  
{m}Y$ makes sense also for $m$ not an integer.  
  
We now discuss the principal symbols of the invariant operators that  
we study. Let  
$$  
        T_{vert}Y:=\ker (TY \to TB)  
$$   
be the bundle of {\em vertical} tangent vectors to $Y$, and let $T_{vert}^*Y$ be  
its dual. We fix compatible metrics on $T_{vert}Y$ and $T_{vert}^*Y$, and  
define $S_{vert}^*Y$, the {\em cosphere bundle of the vertical  
tangent bundle} to $Y$, to be the set of vectors of  
length one of $T_{vert}^*Y$. Also, let   
$$  
        \sigma_m : \Psi^m(Y_b) \to \CI(S_{vert}^* \cap T^*Y_b)  
$$   
be the usual principal symbol map, defined on the space of pseudodifferential   
operators of order $m$ on $Y_b$. The definition of $\sigma_m$ depends on  
the choice of a trivialization of the bundle of homogeneous functions  
of order $m$ on $T^*_{vert}Y$, regarded as a bundle over $S_{vert}^*Y$.  
The principal symbols $\sigma_m(A_b)$ of an element (or family)  
$A=(A_b) \in \Psm {m}Y$ then gives rise to a smooth function on  
$\CI(S_{vert}^*Y)$, which is invariant with respect to $\GR$, and hence  
descends to a smooth function on $S_{vert}^*Y$, which has compact support  
because of the support condition on the kernel of $A$. The resulting   
function,   
\begin{equation}  
        \sigma_m(A) \in \CIc((S_{vert}^*Y)/\GR),  
\end{equation}  
will be referred to as the {\em principal symbol} of an element (or operator)   
in $\Psm {m}Y$.   
  
In the particular case $Y = \GR$, $\Psm {\infty}\GR$ identifies with families of  
convolution operators on the fibers $\GR_b$ with kernels contained in  
a compact subset of $\GR$, smooth outside the identity, and only with  
conormal singularities at the identity. In particular, $\Psm  
{-\infty}{\GR}= \CIc(\GR)$, with the fiberwise convolution product.  
  
Suppose now that the quotient $Y/\GR$ is compact, which implies that  
$(S^*_{vert}Y)/\GR$ is also compact. As it is customary, an operator $A  
\in \Psm {m}Y$ is called {\em elliptic} if, and only if, its principal  
symbol is everywhere invertible. The same definition applies to  
$A = [A_{ij}] \in M_N(\Psm {m}Y):$\   
the operator $A$, regarded as acting on sections  
of the trivial vector bundle $\CC^N$, is elliptic if, and only if,  
its principal symbol    
$$   
        \sigma_m(A) := [\sigma_m(A_{ij})] \in M_n(\CI(S_{vert}^*Y/\GR))   
$$    
is invertible.  
  
Assume that there is given a $\GR$-invariant metric on $T_{vert}Y$,   
the bundle of {vertical} tangent vectors, and a  
$\GR$-equivariant bundle $W$ of modules over the Clifford algebras of  
$T_{vert}Y$. Then a typical example of a family $D=(D_b) \in \Psm {\infty}Y$   
is that of the family of Dirac operators $D_b$ acting on the fibers $Y_b$ of   
$Y \to B$.  (Each $D_b$ acts on sections  
of $W\vert_{Y_b}$, the restriction of the given Clifford module $W$ to  
that fiber.)  
  
Before proceeding, in the next section, to define the equivariant  
family index of an elliptic family invariant with respect to a bundle  
of Lie groups, let us take a closer look at a particular case of the  
previous construction.   
  
Take  
$  
        Y = B \times Z \times \RR^q,  
$  
with $Z$ a compact manifold and   
$  
        \GR = B \times \RR^q,  
$  
$\pi$ and $d$ being the projections onto the first components of  
each product.  The action of $\GR$ on $Y$ is given by translation on  
the last component of $Y$. Then the $\GR$-invariance condition becomes  
simply $\RR^q$ invariance with respect to  
the resulting $\RR^q$ action. If $Y$ and  
$\GR$ are as described here, then we call $Y$ {\em a flat  
$\GR$-space}.    
  
One disadvantage of the algebras $\Psm {\infty}Y$ is the following. It  
is possible to find families $A=(A_b) \in \Psm {0}Y$ such that each   
$A_b$ is invertible as a bounded operator, but the family $(A_b^{-1})$  
is not in $\Psm {0}Y$, although it consists of invariant, pseudodifferential   
operators. This pathology is due to the support condition. Nevertheless,  
for $\GR$ consisting of abelian groups, it is easy to remedy this pathology   
by enlarging the algebra $\Psm {-\infty}Y$, as follows.     
   
Since the enlargement of the algebra $\Psm {-\infty} Y$ is done  
locally, we may assume that $Y$ is a flat $\GR$--space.  The residual  
ideal of the algebra $\Psm {\infty}Y$ is $\Psm {-\infty}Y$ and  
consists of operators that are regularizing along each fiber. More  
precisely, it consists of those families of smoothing operators on $Y  
= B \times Z \times\RR^q$ that are translation-invariant under the  
action of $\RR^q$ and have {\em compactly supported} convolution  
kernels. Thus  
$$  
        \Psm {-\infty}Y \cong \CIc(B \times \RR^q; \Psi^{-\infty}(Z))  
        \subset \mathcal{S}(B \times \RR^q; \Psi^{-\infty}(Z)) \simeq  
        \mathcal{S}(B \times Z\times Z\times\RR^q).  
$$  
(Here ${\mathcal S}$ is the generic notation for the space of Schwartz  
functions on a suitable space, in this case on $B \times \RR^q$, with  
values regularizing operators.) The second isomorphism above is  
obtained from the isomorphism  
$$  
        \Psi^{-\infty}(Z) \simeq \CI(Z \times Z),  
$$  
defined by the choice of a nowhere vanishing density on $Z$. If we  
also endow $\Psi^{-\infty}(Z)$ with the locally convex topology  
induced by this isomorphism, then it becomes a nuclear, locally convex  
space.  We now enlarge the algebra $\Psm {\infty}Y$ to include all  
invariant, regularizing operators whose kernels are in $\mathcal{S}(B  
\times Z\times Z\times\RR^q)$:  
\begin{equation}  
        \PsS{m}Y = \Psm{m}Y + \mathcal{S}(B \times Z\times  
        Z\times\RR^q).  
\end{equation}   
Explicitly, the action of $T \in \mathcal{S}(B \times Z\times  
Z\times\RR^q)$ on a smooth function $f \in \CIc(B \times Z \times  
\RR^q)$ is given by  
$$  
        Tf(b,y_1,t) = \int_{Z \times \RR^q} T(b, y_1,  
        y_0,t-s)f(b,y_0,s) dy_0 ds.  
$$  
For $B$ reduced to a point, the algebras $\PsS {\infty}Y$ were  
introduced in \cite{Melrose42} as the range space of the indicial map  
for pseudodifferential operators on manifolds with corners  
  
Still assuming that we are in the case of a flat $\GR$-space, we  
notice that the Fourier transformation ${\mathcal F}_2$ in the  
translation-invariant directions gives a dual identification  
\begin{equation}  
        \PsS {-\infty}Y \ni T \to \hat{T}:=  
        {\mathcal F}_2 T {\mathcal F}_2^{-1}  
        \in \mathcal{S}(B \times Z\times Z \times \RR^q) =  
        \mathcal{S}(B \times \RR^q; \Psi^{-\infty}(Z))  
\label{Smoothing.Ideal}  
\end{equation}  
with the space of Schwartz functions with values in the smoothing  
ideal of $Z.$ The point of this identification is that the  
convolution product is transformed into the pointwise product. The  
Schwartz topology from \eqref{Smoothing.Ideal} then gives $\PsS  
{-\infty}Y$ the structure of a nuclear locally convex topological  
algebra. Following the same recipe, the Fourier transform also gives  
rise to an {\em indicial family}  
\begin{equation}\label{Fourier.Transform}  
        \Phi:\PsS {\infty}Y\ni T \to \hat{T}:={\mathcal F}_2 T {\mathcal F}_2^{-1}  
        \in \CI(B \times \RR^q;\Ps{\infty}{Z}).  
\end{equation}  
We denote   
$$  
        \hat{T}(\tau)=\Phi(T)(\tau).  
$$  
The map $\Phi$ is not an isomorphism since $\hat{A}(\tau)$ has joint  
symbolic properties in the variables of $\RR^q$ and the fiber  
variables of $T^*Z.$ Actually, the principal symbols of the  
operators $\hat{A}(\tau)$ is constant on the fibers of $T_{vert}^*Y \to  
B$.

\begin{lemma}\label{Lemma.dilations}\ Assume $Y=B \times Z \times \RR^q$   
is a flat $\GR$-space. Then the action of the group $GL_q(\RR)$ on the  
last factor of $Y = B \times Z \times \RR^q$ extends to an action by  
automorphisms of $\CI(B,GL_q(\RR))$ on $\Psm {\infty} Y$ and on $\PsS  
{\infty} Y.$  
\end{lemma}

\begin{proof}\   
The vector representation of $GL_q(\RR)$ on the second component of  
$Z \times \RR^q$ defines an action of $GL_q(\RR)$ on  
$\Psi^\infty(Z \times \RR^q)$ that preserves the class of properly  
supported operators and the products of such operators.  It also  
normalizes the group $\RR^q$ of translations, and hence it maps  
$\RR^q$-invariant operators to $\RR^q$-invariant operators. This  
property extends right away to the action of $\CI(B,GL_q(\RR))$ on  
families of operators on $B \times Z \times \RR^q$, and hence  
$\CI(B,GL_q(\RR))$ maps $\Psm {m}Y$ isomorphically to itself. From the  
isomorphism \eqref{Smoothing.Ideal}, we see that $\CI(B,GL_q(\RR))$  
also maps $\PsS {-\infty}Y$ to itself. This gives an action by  
automorphisms of $\CI(B,GL_q(\RR))$ on $\PsS {\infty} Y$, which is the  
sum of $\PsS {-\infty}Y$ and $\Psm {\infty}Y$.  
\end{proof}

Suppose the family of Lie groups $\GR$ consists of abelian Lie groups,  
so that $\GR$ is a vector bundle, $\GR \cong \lgg$. By choosing a lift  
of $Y /\GR \to Y$, which is possible because the fibers are  
contractable, we obtain that locally the bundle $Y$ is isomorphic to a  
flat $\GR$ space.  Then the above lemma allows us to extend the  
previous definitions, including those of the algebras $\PsS {\infty}Y$  
and of the indicial family from the flat case to the case $\GR$  
abelian. The indicial family $\hat{A}$ of an operator $A \in \PsS  
{\infty}Y$, will then be a family of pseudodifferential operators  
acting on the fibers of $Y \times_B \lgg^* \to \lgg^*$ (here $\lgg^*$  
is the dual of the vector bundle $\lgg$):  
\begin{equation}\label{eq.bundle.hatA}   
        \hat{A}(\tau) \in \Psi^*(Y_b/\GR_b), \quad \text{if } \tau \in  
        \lgg_b^*.  
\end{equation}   
The action of $GL_q(\RR)$ in the above lemma will have to be replaced with   
the group of gauge-transformations of $\lgg$.  
  
We can look at the general abelian $\GR$ from the point an other, related  
point of view also. Let $Z := Y/\GR$, which is again a fiber bundle  
$Z \to B$. The algebra $\Psm{\infty}Z$ of smooth families of pseudodifferential   
operators along the fibers of $Z \to B$ can be regarded as the algebra of  
sections of a vector bundle on $Z$ (with infinite dimensional fibers).  
We can lift this vector bundle to $\lgg^*$ and then $\hat{A}$ is a section   
of this lifted vector bundle over $\lgg^*$, and we write this as follows:   
\begin{equation}\label{eq.section.hatA}   
	\hat{A} \in \Gamma(\lgg^*,\Psm{\infty}Z).  
\end{equation}  
   
The considerations of this section extend immediately to operators  
acting between sections of a $\GR$-equivariant vector bundles. If  
$E_0$ and $E_1$ are $\GR$-equivariant vector bundles, we denote by  
$\psi^m(Y;E_0,E_1)$ the space of $\GR$-invariant families of order $m$  
pseudodifferential operators acting on sections of $E_0$, with values  
sections of $E_1$, whose convolution kernels have compact support.

\section{The equivariant family index: Definition\label{S.Homotopy.Inv}}

We now define study invariants of elliptic operators in $M_N(\Psm  
{\infty}Y)$, the main invariant being the (equivariant family) index  
of such an invariant, elliptic family. For $\GR$ consisting of simply  
connected, solvable Lie groups and $\dim Y > \dim \GR$, we then show  
that the (equivariant family) index gives the obstruction for family  
$A=(A_b) \in M_N(\Psm{m}Y)$ to have a perturbation by a family  
$R=(R_b) \in M_N(\Psm{-\infty}Y)$, such that $A+R = (A_b + R_b)$ be  
invertible, for all $b \in B$, between suitable Sobolev spaces, see  
Theorem \ref{Theorem.obst}. For families of abelian Lie groups $\GR$,  
we give an interpretation of the index of an elliptic operator in  
terms of its indicial family. This leads to an Atiyah-Singer index  
type formula for the Chern character of the index of a family of  
invariant, elliptic operators. If $\GR$ is abelian (that is, if its  
fibers are abelian Lie groups), then we can consider the algebra $\PsS  
{\infty}Y$ instead of $\Psm {\infty}Y$.  
  
We now proceed to define the index of an elliptic family $A \in \Psm {m}Y$.  
This will be done using the $K$-theory of Banach algebras. Let $C^*_r(Y,\GR)$  
be the closure of $\Psm {-\infty}{Y}$ with respect to the norm  
$$  
        \| A \| = \sup_{b \in B} \|A_b\|  
$$  
each operator $A_b$ acting on the Hilbert space of square integrable  
densities on the fiber $Y_b$. If $Y = \GR$, then we also write  
$C^*_r(\GR,\GR)=C^*_r(\GR)$, $C^*_r(\GR)$ is the reduced $C^*$-algebra  
associated to $\GR$, regarded as a groupoid. For each locally compact  
space $X$, we denote by $C_0(X)$ the space of continuous functions  
vanishing at infinity on $X$.  Then, if $\GR$ is abelian, we have  
$C^*_r(\GR) \simeq C_0(\GR^*)$.  
  
We shall use below $\widehat{\otimes}$, the minimal tensor product of  
$C^*$-algebras. This minimal tensor product is defined to be  
(isomorphic to) the completion of the image of $\pi_1 \otimes \pi_2$,  
the tensor product of two injective representations $\pi_1$ and  
$\pi_2$.  For the cases we are interested in, the minimal and the  
maximal tensor product coincide \cite{Sakai}.

\begin{lemma}\label{lemma.stable}\   
Assume $\dim Y > \dim \GR$. Also, let ${\mathcal K}={\mathcal  
K}(Y_b/\GR_b)$ denote the algebra of compact operators on one of the  
fibers $Y_b/\GR_b$, for some fixed but arbitrary $b \in B$. Then  
$$  
        C_r^*(Y,\GR) \simeq C_r^*(\GR) \widehat{\otimes} {\mathcal K}.  
$$  
Consequently, $K_i(C_r^*(Y,\GR)) \simeq K_i(C_r^*(\GR))$.  
\end{lemma}

\begin{proof}\ Let $\mfk A$ be the space of sections of  
the bundle of algebras ${\mathcal K}(Y_b/\GR_b)$.  If $Y$ is a flat  
$\GR$-space, then the isomorphism $C_r^*(Y,\GR) \simeq C_r^*(\GR)  
\widehat{\otimes} {\mathcal K}$ follows, for example, from the results  
of \cite{LauterNistor}. In general, this local isomorphism gives 
$C_r^*(Y,\GR) \simeq C_r^*(\GR) \widehat{\otimes}_{C_0(B)} {\mathfrak 
A}$. 
  
Our assumptions imply that ${\mathcal K}$ is infinite dimensional, and  
hence its group of automorphisms is contractable, see  
\cite{Diximier}. Consequently, there is no obstruction to trivialize  
the bundle of algebras ${\mathcal K}(Y_b/\GR_b)$, and hence  
${\mathfrak A} \simeq C_0(B) \widehat{\otimes}{\mathcal K}$. We obtain  
\begin{equation*}   
        C_r^*(Y,\GR) \simeq C_r^*(\GR) \widehat{\otimes}_{C_0(B)} 
        {\mathfrak A} \simeq C_r^*(\GR) \widehat{\otimes}_{C_0(B)} 
        C_0(B) \widehat{\otimes}{\mathcal K} \simeq C_r^*(\GR) 
        \widehat{\otimes} {\mathcal K}. 
\end{equation*}    
   
The last part of lemma follows from the above results and from the  
natural isomorphism $K_i(A \widehat{\otimes} {\mathcal K}) \simeq  
K_i(A)$, valid for any $C^*$-algebra $A$.  
\end{proof}

We proceed now to define the index of an {\em elliptic, invariant  
family} of operators  
$$   
        A = (A_b) \in M_N(\Psm{m}Y) = \Psm{m}{Y;\CC^N}.  
$$   
We assume first that $Y/\GR$ is compact, for simplicity; otherwise, we  
need to use algebras with adjoint units. We observe that there exists  
an exact sequence  
\begin{equation}\label{eq.exact.seq}  
  	0 \to C^*_r(Y,\GR) \to {\mathcal E} \to \CI(S^*_{vert}Y) \to  
  	0,\quad {\mathcal E}:= \Psm {0}Y + C^*_r(Y,\GR),  
\end{equation}  
obtained using the results of \cite{LauterNistor}. The operator $A$  
(or, rather, the family of operators $A = (A_b)$) has an invertible  
principal symbol, and hence the family $T = (T_b)$,  
$$   
        T_b:=(1 + A_b^*A_b)^{-1/2}A_b,  
$$  
consists of elliptic, invariant operators. Moreover, $T$ is an element  
of ${\mathcal E} := \Psm {0}Y + C^*_r(Y,\GR)$, its principal symbol is  
still invertible, and hence defines a class $[T] \in  
K_1(\CI(S^*_{vert}Y)) \simeq K^1(S^*_{vert}Y)$. Let  
\begin{equation*}\label{eq.boundary.map}  
        \pa : K_1^{alg}(S^*_{vert}Y) \to K_0^{alg}(C^*_r(Y,\GR))  
        \simeq K_0 (C^*_r(Y,\GR))  
\end{equation*}  
be the boundary map in the $K$-theory exact sequence  
\begin{multline*}  
  	K_1^{alg}(C^*_r(Y,\GR)) \to K_1^{alg}({\mathcal E}) \to  
  	K_1^{alg}(S^*_{vert}Y)  
  	\stackrel{\pa}{\rightarrow}K_0^{alg}(C^*_r(Y,\GR)) \\ \to  
  	K_0^{alg}({\mathcal E}) \to K_0^{alg}(S^*_{vert}Y)  
\end{multline*}  
associated to the exact sequence \eqref{eq.exact.seq}. Because $K_0 
(C^*_r(Y,\GR)) \simeq K_0(C^*_r(\GR)),$ by Lemma \ref{lemma.stable}, 
we get a group morphism 
\begin{equation}  
        \ind_a : K_1^{\alg}(\CI(S^*_{vert}Y))  \to  K_0(C^*_r(\GR)),  
\end{equation}  
which we shall call {\em the analytic index morphism}. The image of  
$A$ under the composition of the above maps is $\ind_a([T])$, will be  
denoted $\ind_a(A)$, and will be called {\em the analytic index of}  
$A$. A more direct but longer definition is contained in the proof of  
Theorem \ref{Theorem.obst} below (see Equation \eqref{eq.def.ind}).  
The analytic index morphism descends in this case to a group morphism  
$K_1^{\alg}(\CI(S^*_{vert}Y)) \to K_0(C^*_r(\GR))$ denoted in the same  
way.  For $\GR$ abelian, we can replace ${\mathcal E}$ with $\PsS  
{0}Y$ and $C^*_r(\GR)$ with $\PsS {-1}{Y}$.  
  
We denote by $I_N$ the unit of the matrix algebra $M_N(\mathcal E)$.  
If $A \in \Psm{0}{Y;E}$, we can find $N$ large such that $\Psm{0}{Y;E}  
\subset M_N(\Psm{0}Y)$ with $1$ mapping to the projector $e$ under  
this isomorphism. Then $A + (I_N-e)$ is invertible in this matrix  
algebra, and we define then $\ind_a(A) = \ind_a(A + I_N - e)$.  If  
$Y/\GR$ is not compact, this definition of the index applies only to  
elliptic operators in $I_N + M_N(\Psm{0}{Y;E})$. All results below  
extend to these operators, after obvious changes.  
  
For differential operators acting between sections of different  
bundles, we can define the analytic index using the adiabatic groupoid  
of $\GR$ as in \cite{LauterNistor}. For arbitrary elliptic  
pseudodifferential operators $A \in \Psm{m}{Y;E_0,E_1}$, acting  
between sections of possibly different vector bundles, we can define  
the index using Kasparov's bivariant $K$-Theory. Since for $Y/\GR$  
non-compact, no element in $\Psm{m}{Y;E_0,E_1}$ is invertible modulo  
regularizing operators, we must allow more general operators in this  
case. For example, we can take $A \in \Gamma(Y,\Hom(E_0,E_1))^{\GR} +  
\Psm{0}{Y;E_0,E_1}$, provided that it is bounded. If $A \in  
\Psm{m}{Y;E_0,E_1}$, $m >0$, we replace $A$ by $(1 + AA^*)^{-1/2} A$,  
which will be an operator in the closure of $\Psm{0}{Y;E_0,E_1}$, by  
the results of \cite{LauterNistor2}. If $A$ is elliptic, then the  
resulting family defines by definition an element of $KK(\CC,  
C^*(\GR))$.  To obtain the index, we use the isomorphism $KK(\CC,  
C^*(\GR)) \cong K_0(C^*(\GR))$. This definition has the disadvantage  
that we must always work with $C^*(\GR)$, the closure of  
$\Psm{0}{\GR}$, and we cannot use $\PsS{-\infty}{\GR}$.  
  
The main property of the analytic index of an operator $A$ is that it 
gives the obstruction to the existence of invertible perturbations of 
$A$ by lower order operators. We denote by $H^s(Y_b)$ the $s$th 
Sobolev space of $1/2$-densities on $Y_b$, which is uniquely defined 
because of the bounded geometry of $Y_b$, for $Y/\GR$ compact. More 
precisely, $H^s(Y_b)$ is, by definition, the domain of $(1 + 
D^*D)^{s/2m}$, if $D \in 
\Psm{m}Y$ is elliptic and $s \ge 0$. For $s \le 0$, $H^s(Y_b)$ is, by  
definition, the dual of $H^{-s}(Y_b)$.

\begin{theorem}\label{Theorem.obst}\   
Let $\GR \to B$ be a bundle of Lie groups acting on the fiber bundle  
$Y \to B$, as above, and assume that $Y/\GR$ is compact, of positive  
dimension.  Let $A \in \Psm {m}{Y,\CC^N}$ be an elliptic operator.  
Then we can find $R \in \Psm {m-1}{Y,\CC^N}$ such that  
$$   
        A_b + R_b : H^{s}(Y_b)^N \to H^{s-m}(Y_b)^N   
$$    
is invertible for all $b \in B$ if, and only if, $\ind_a(A)=0$.  
Moreover, if $\ind_a(A)=0$, then we can choose $R \in \Psm{-\infty}Y$.  
The same result holds if $Y/\GR$ is non-compact and  
$A \in \Gamma(Y,\Hom(E_0,E_1))^{\GR} + \Psm {0} {Y;E_0,E_1}$ is elliptic  
and bounded.   
\end{theorem}

\begin{proof}\   
It is clear from definition that if we can find $R$ with the desired  
properties, then $\ind_a(A) = 0 \in K^0(\lgg^*)$.  Suppose now that $A  
\in \Psm {m}{Y,\CC^N}$ is elliptic and has vanishing analytic index.  
Using the notation $T = (1 + AA^*)^{-1/2}A$, we see that $A_b$ is  
invertible between the indicated Sobolev spaces if, and only if, $T_b$  
is invertible as a bounded operator on $L^2(Y_b)^N$. We can hence  
assume that $m = 0$ and $T = A$.  
  
Because $C^*_r(Y,\GR)$ satisfies $C^*_r(Y,\GR)\simeq C^*_r(\GR)  
\widehat{\otimes} {\mathcal K}$, by Lemma \ref{lemma.stable}, we can  
use some general techniques to prove that the vanishing of $\ind_a(A)$  
implies that $A$ has a perturbation by invariant, regularizing  
operators in $\Psm{-\infty}{Y,\CC^N}$ that is invertible on each  
fiber. We fix an isomorphism $M_N(C^*_r(Y,\GR)) \cong C^*_r(\GR)  
\otimes {\mathcal K}$. We now review this general technique using a  
generalization of an argument from \cite{NistorKthry}. Let $\mathcal  
E$ be the algebra introduced in Equation \eqref{eq.exact.seq}.  We  
denote by $I_N$ the unit of the matrix algebra $M_N(\mathcal E)$.  
Also, denote by $\overline{\mathcal E}$ the closure of $\mathcal E$ in  
norm.  
  
Choose a sequence of projections $p_n \in {\mathcal K}$, $\dim p_n =  
n$, such that $p_n \to 1$ in the strong topology. Because $A_b$ is  
invertible modulo $C^*_r(\GR_b) \otimes {\mathcal K}$, we can find a  
large $n$ and $R \in M_N(\Psm{-\infty}Y)$ such that  
$$  
	A'_b := A_b \oplus R_b : L^2(Y_b)^N \oplus L^2(Y_b)^N \to L^2(Y_b)^N  
$$   
is surjective and $(1 \otimes p_n) R_b (1 \otimes p_n) = R_b$, for all  
$b \in B$. Then $\{0\}$ is not in the spectrum of ${A'}{A'}^*$, and we  
can consider $V := (A'{A'}^*)^{-1/2} A' \in M_N(\overline{\mathcal  
E})$, which by construction satisfies $VV^* = I_N \in  
M_N(\overline{\mathcal E})$. Consequently, $V^*V$ is a projection in  
$M_{2N}(\overline{\mathcal E})$. Because $A$ is invertible modulo  
$M_N(C^*_r(\GR) \otimes {\mathcal K})$,  
$$  
	V^*V - VV^* \in  M_{2N}(C^*_r(\GR) \otimes {\mathcal K}).  
$$  
Let $e = I_N \oplus (1 \otimes p_n) -V^*V$, which is also a  
projection, by construction. Moreover,  
\begin{equation}\label{eq.def.K}  
	e - (1 \otimes p_n) \in M_{2N}(C^*_r(\GR) \otimes \mathcal K),  
\end{equation}  
and hence both $e$ and $1 \otimes p_n$ are projection in   
$M_{2N}(C^*_r(\GR)^+ \otimes \mathcal K)$ (for any algebra $B$,   
we denote by $B^+$ the algebra with adjoint unit). Equation \eqref{eq.def.K}  
gives that, by definition, $[e] - [1 \otimes p_n]$ defines a $K$-theory class  
in $K_0(C^*_r(\GR))$. From definition, we get then  
\begin{equation}\label{eq.def.ind}  
	\ind_a(A) = [e] - [1 \otimes p_n].  
\end{equation}  
  
Now, if $\ind_a(A)=0$, then we can find $k$ such that $e \oplus  
1\otimes p_k$ is (Murray-von Neumann) equivalent to $1\otimes p_{n +  
k}$. By replacing our original choice for $n$ with $n + k$, we may  
assume that $e$ and $1 \otimes p_n$ are equivalent, and hence that we  
can find $u \in \CC I_{2N} + M_{2N}(C^*_r(\GR) \otimes {\mathcal K})$  
with the following properties:\ there exists a large $l$ and $x \in  
M_{2N}(C^*_r(\GR))$ such that, if we denote $e_0 = 1 \otimes p_l  
\oplus 1 \otimes p_n$, then $u = I_{2N} + e_0 x e_0$ and $e = u(1  
\otimes p_n) u^{-1}$. Then $Vu$ is in $M_n (\overline{\mathcal E})$  
(more precisely $I_N V I_N = V$) and is invertible.  Consequently  
$B:=(A'{A'}^*)^{1/2} V$ is also invertible. But $B$ is a perturbation  
of $A'$, and hence also of $A$, by an element in $M_{2N}(C^*_r(\GR)  
\otimes {\mathcal K})$. Since $\Psm{-\infty}Y$ is dense in  
$C^*_r(\GR)$, this gives the result.  
\end{proof}

\section{The index for bundles of solvable Lie groups\label{Sec.Solvable}}

We now treat in more detail the case of bundles of solvable Lie  
groups, when more precise results can be obtained. Other classes of  
groups will lead to completely different problems and results, so we  
leave their study for the future. The class of simply-connected  
solvable fibers is rich enough and has many specific features, so we  
content ourselves from now on with this case only.  \\[2mm]  
{\bf Assumption.}\ {\em From now on and throughout this paper, we  
shall assume that the family $\GR$ consists of simply-connected  
solvable Lie groups.} By ``simply-connected'' we mean, as usual,  
``connected with trivial fundamental group.''  
  
\vspace*{2mm} We shall denote by $\lgg \to B$ the bundle of Lie  
algebras of the groups $\GR_b$, $\lgg_b \simeq Lie (\GR_b)$ and by  
$$  
        \exp: \lgg \to \GR.  
$$  
the exponential map.   
  
Because all the groups $\GR_b$ are solvable, we have that the enveloping  
$C^*$ algebra of $\GR$, that is, $C^*(\GR)$ is isomorphic to the reduced  
$C^*$-algebra of $\GR$: \ $C^*(\GR) \simeq C^*_r(\GR)$, so we can drop the  
index ``$r$'' from the notation.  
  
In order to study the algebra $C^*_r(\GR) = C^*(\GR) $ and its $K$-groups, we  
shall deform it to a commutative algebra. This deformation is obtained  
as follows.  Let $\GR_{ad} = \{0\} \times \lgg \cup (0,1] \times \GR$,  
$B_1 = [0,1] \times B$, and $d:\GR_{ad} \to B_1$ be the natural  
projection. On $\GR_{ad}$ we put the smooth structure induced by  
$$  
        \phi: B_1 \times \lgg \to \GR_{ad}  
$$  
$\phi(0,X) = (0,X)$ and $\phi(t,X)=(t,\exp(tY))$, which is a bijection  
for all $(t,X) \in [0,1] \times \lgg$ in a small neighborhood of $\{0  
\} \times \lgg \cup B_1 \times \{0\}$.  Then we endow $\GR_{ad}$ with  
the Lie bundle structure induced by the pointwise product.  Evaluation  
at $t \in [0,1]$ induces algebra morphisms  
$$  
        e_t : \Psm {m}{\GR_{ad}} \to \Psm {m}{\GR}, \quad t > 0,  
$$  
and   
$$  
        e_0 : \Psm {m}{\GR_{ad}} \to \Psm {m}{\lgg}, \quad t = 0.  
$$  
Passing to completions, we obtain morphisms $e_t$ from $C^*(\GR_{ad})$  
to $C^*(\GR)$, for $t > 0$, and to $C^*(d^{-1}(0) \times B) \simeq  
C_0(\lgg^*)$, for $t = 0$. (This construction should be compared to that  
of the normal groupoid \cite{Skandalis}.)

\begin{lemma}\label{Lemma.Hom.Inv}\   
The morphisms $e_t : C^*(\GR_{ad}) \to C^*(\GR)$, for $t > 0$, and  
$e_0 : C^*(\GR_{ad}) \to C_0(\lgg^*)$, for $t = 0$, induce  
isomorphisms in $K$-theory.  
\end{lemma}

\begin{proof}\ Assume first that there exists a Lie group bundle morphism   
$\GR \to B \times \RR$. (In other words, there exists a smooth map  
$\GR \to \RR$ that is a morphism on each fiber.) Let $\GR'$ denote the  
kernel of this morphism and let $\GR_{ad}'$ be obtained from $\GR'$ by  
the same deformation construction by which $\GR_{ad}$ was obtained  
from $\GR$. Then we obtain a smooth map $\GR_{ad} \to \RR$ that is a  
group morphism on each fiber, and hence  
$$   
        C^*(\GR_{ad}) \simeq C^*(\GR_{ad}') \rtimes \RR, \; C^*(\GR)  
        \simeq C^*(\GR') \rtimes \RR, \; \text{ and } C_0(\lgg^*)  
        \simeq C_0({\lgg'}^*) \rtimes \RR.  
$$   
Moreover, all above isomorphisms are natural, and hence compatible with   
the morphisms $e_t$. Assuming now that the result was proved for all   
Lie group bundles of smaller dimension, we obtain the desired result   
for $\GR$ using Connes' Thom isomorphism in $K$-theory \cite{ConnesThom},    
which in this particular case gives:    
$$   
        K_i(C^*(\GR_{ad}))  \simeq K_{i+1} (C^*(\GR_{ad}')), \;    
        K_i(C^*(\GR)) \simeq K_{i+1}(C^*(\GR')), \; \text{ and }  
$$  
$$   
        K_i(C_0(\lgg^*)) \simeq K_{i+1} (C_0({\lgg'}^*)).   
$$   
   
This will allow us to complete the result in the following way. Let  
$U_k$ be the open subset of $B$ consisting of those $b \in B$ such  
that $[\GR_b,\GR_b]$ has dimension $\ge k$. From the Five Lemma and  
the six term exact sequences in $K$-theory associated to the ideal  
$C^*(\GR_{ad}\vert_{U_{k} \smallsetminus U_{k+1}})$ of  
$C^*(\GR_{ad}\vert_{U_{k-1} \smallsetminus U_{k+1}})$, for each $k$, we see  
that it is enough to prove our result for $\GR_{ad}\vert_{U_k  
\smallsetminus U_{k+1}}$ for all $k$. Thus, by replacing $B$ with $U_k  
\smallsetminus U_{k+1}$, we may assume that the rank of the abelianization  
of $\GR_b$ is independent of $b$.  Consequently, the abelianizations  
of $\GR_b$ form a vector bundle  
$$   
        {\mathcal A} := \cup \GR_b / [\GR_b,\GR_b]   
$$   
on $B$.   
   
A similar argument, using the Meyer-Vietoris exact sequence in  
$K$-theory and the compatibility of $K$-theory with inductive limits  
\cite{Blackadar}, shows that we may also assume the vector bundle  
${\mathcal A}$ of abelianizations to be trivial. Then the argument at  
the beginning of the proof applies, and the result is proved.  
\end{proof}

{}From the above lemma we immediately obtain the following corollary:

\begin{corollary}\label{Cor.Hom.Inv}\   
Let $\GR$ be a bundle of simply connected, solvable Lie groups. Then  
$$  
        K_i(C^*(\GR)) \simeq K_i(C_0(\lgg^*)) \simeq K^i(\lgg^*).  
$$  
\end{corollary}

We now give an interpretation of $\ind_a(A)$, for $\GR$ abelian, using  
the properties of the indicial family $\hat{A}(\tau)$ of $A$. We assume  
that $Y/\GR$ is compact.  
   
We shall also use the following construction.  Let $X$ be a compact  
manifold with boundary.  Let $T(x)$ be a family of elliptic  
pseudodifferential operators acting between sections of two smooth  
vector bundles, $E_0$ and $E_1$, on the fibers of a fiber bundle $M  
\to X$ whose fibers are compact manifolds without corners. Then we can  
realize the index of $T$ as an element in the relative group  
$K^0(X,\pa X)$. This can be done directly using Kasparov's theory or  
by the ``Atiyah-Singer trick'' as follows. We  and proceed as in  
\cite{AS4}, Proposition (2.2), to define a smooth family of maps  
$R(x): \CC^N \to \CI(Y)$, such that the induced map  
$$  
	V := T \oplus R : \CI(X)^N \oplus \CI(X,L^2(Y;E_0)) \to   
	\CI(X,L^2(Y;E_1))  
$$  
is onto for each $x$. Since $T(x)$ is invertible for $x \in \pa X$, we  
can choose $R(x)=0$ for $x \in \pa X$. Then $\ker (V)$ is a vector  
bundle on $X$, which is canonically trivial on the boundary $\pa X$.  
The general definition of the index of the family $T$ in \cite{AS4} is  
that of the difference of the kernel bundle $\ker (V)$ and the trivial  
bundle $X\times \CC^N$. Since the bundle $\ker (V)$ is canonically  
trivial on the boundary of $X$, we obtain an element  
$$  
	\deg(T) \in K^0(X,\pa X).  
$$  
The degree is invariant with respect to homotopies $T_t$ of families  
of operators on $X$ that are {\em invertible} on $\pa X$ throughout  
the homotopy. We shall use the degree in Theorem \ref{Theorem.degree}  
for $X=B_R$, a large closed ball in $\GR^* \simeq \lgg^*$, or for $X$  
being the radial compactification of $\lgg^*$. If the boundary of $X$  
is empty, this construction goes back to Atiyah and Singer and then  
$\deg(T)$ is simply the family index of $T$. If $\pa X$ is not empty,  
this definition of $\deg(T)$ is due to Melrose.  We note that when  
$\pa X \not = \emptyset$, the degree is not a local quantity in $T$,  
in the sense that it depends on more than just the principal symbol.  
  
Assume now that the family $T$ above consists of order zero operators  
and $T(x)$ is a multiplication operator for each $x \in \pa X$. We  
want to compute the Chern character of $\deg(T)$ using the  
Atiyah-Singer family index formula \cite{AS1,AS4}.  To introduce the  
main ingredients of the index formula, denote by $S^*_{vert}M$ the set  
of unit vectors in the dual of the vertical tangent bundle $T_{vert}M$  
to the fibers of $M \to X$.  Because the family $T$ is elliptic, the  
principal symbols define an invertible matrix of functions  
$$  
	a = \sigma_0(T) \in \CI(S^*_{vert}M;\Hom(E_0,E_1)).  
$$   
Since the operators $T(x)$ are multiplication operators, we can then  
extend $a$ to an invertible endomorphism on  
$S_M:=S^*_{vert}M \cup B^*$, with $B^*$ denoting the set of vertical  
cotangent vectors of length $\le 1$ above $\pa X$, as in  
\cite{AB}.   
The constructions of \cite{AS1,AS4} are in terms of    
$[a']\in K^0(T^*_{vert}M, T^*_{vert} M\vert_{\pa M})$ obtained  
by applying the clutching (or difference) construction to  
$a$. Explicitly, $[a']$ is represented by $(E_0,E_1,a_1)$ (where $a_1$ is  
a smooth function that coincides with $a$ outside a neighborhood of the  
zero section). It defines an element in  
$$  
	K^0(T^*_{vert}M, T^*_{vert} M\vert_{\pa M}) = K^0(T^*_{vert}M  
	\smallsetminus T^*_{vert} M\vert_{\pa M}),  
$$   
because $a$ defines an endomorphism of the trivial bundle $\CC^N$  
which is invertible outside a compact set (see \cite{AtiyahKTHRY}).   
  
When $E_0 \cong E_1$, we can assume that  
$E_0 = E_1 = M \times \CC^N$ are trivial of rank $N$, and we have  
that $a$ is an invertible matrix valued function, which hence  
defines an element $[a] \in K^1(S_M)$.    
Let $B_1$ be the set of vectors of norm at most  
$1$ in $T^*M$. After the identification (up to homeomorphism) of $B_1  
\smallsetminus S_M$, the interior of $S_M$, with the difference  
$T^*_{vert}M \smallsetminus T^*_{vert} M\vert_{\pa M}$, we have  
$$  
	[a'] = \pa [a].  
$$  
  
We denote by $\pi_* : H^*(S_M) \to H^{* - 2n + 1}(X, \pa X)$ the  
integration along the fibers, where $n$ is the dimension of the fibers  
of $M \to X$. Integration along the fibers in this case is the  
composition of  
$$  
	\pa : H^*(S_M) \to H^{* + 1} (B_1, S_M) \cong  
	H_c^{*+1}(T^*_{vert}M \smallsetminus T^*_{vert} M\vert_{\pa M})  
$$   
and   
$$  
	\tilde\pi_*: H_c^*(T^*_{vert}M \smallsetminus T^*_{vert}  
	M\vert_{\pa M}) \to H_c^{* -2n}(X \smallsetminus \pa X)  
$$  
obtained by integration along the fibers of the bundle $T^*_{vert}M  
\smallsetminus T^*_{vert} M\vert_{\pa M} \to X \smallsetminus \pa X$:  
$$  
	\pi_* = \tilde \pi_* \circ \pa.  
$$  
  
To state the following result, we also need $Ch: K_1(S_M) \to  
H^{odd}(S_M)$, the Chern character in $K$-Theory and $\mathcal T$, the  
Todd class of $(T^*_{vert}M) \otimes \CC$, the complexified vertical  
tangent bundle of the fibration $M \to X$, as in \cite{AS4}.  Using  
the notation introduced above, we have:

\begin{theorem}\label{Theorem.mult.b}\ Let $M \to X$ be a smooth  
fiber bundle whose fibers are smooth manifolds (without corners), and  
let $T$ be a family of order zero elliptic pseudodifferential  
operators acting along the fibers of $M \to X$. Assume $X$ is a  
manifold with boundary $\pa X$ such that the operators $T(x)$ are  
multiplication operators on $\pa X$, also let $[a']$ and $[a]$ be the  
classes defined above. Then  
\begin{equation*}  
	Ch(\deg(T)) = (-1)^n \tilde\pi_*\big (Ch[a'] {\mathcal  
	T}\big) \in H^{*}(X, \pa X),  
\end{equation*}  
If $E_0 \cong E_1$, then we also have $Ch(\deg(T))= (-1)^n \pi_*\big  
(Ch[a] {\mathcal T}\big)$.  
\end{theorem}

\begin{proof}\ For continuous families $T(x)$ that are multiplication   
operators on the boundary $\pa X$, the degree is a local quantity --  
it depends only on the principal symbol -- so we can follow word for  
word \cite{AS4} to prove that  
$$  
        Ch(\deg(T)) = (-1)^n \tilde\pi_*\big (Ch[a'] {\mathcal  
        T}\big) \in H^{*}(X, \pa X).  
$$  
  
When $E_0 \cong E_1$, using $Ch [a'] = Ch (\pa [a]) = \pa Ch [a]$, we get  
\begin{multline*}  
	Ch(\deg(T)) = (-1)^n \tilde\pi_* \big ( \pa Ch[a] {\mathcal  
      	T} \big )  = (-1)^n \tilde \pi_* \circ \pa (Ch[a]  
      	{\mathcal T}) \\  = (-1)^n \pi_*\big (Ch[a] {\mathcal  
      	T}\big) \in H^{*}(X, \pa X).  
\end{multline*}  
This completes the proof.  
\end{proof}

It is interesting to note that it is not possible in general to give a  
formula for $Ch(\deg(T))$ only in terms of its principal symbol,  
without further assumptions on $T$. A consequence is that, in general,  
the formula for $Ch(\deg(T))$ will involve some non-local  
invariants. It would be nevertheless useful to find such a formula.  
   
Returning to our considerations, we continue to assume that $Z:=Y  
/\GR$ is compact, and we fix a metric on $\GR$ (which, we recall, is a  
vector bundle in these considerations). If $A \in \PsS{\infty}{Y;E_0,E_1}$   
is elliptic (in the sense that its principal symbol is invertible outside  
the zero section), then the indicial operators $\hat{A}(\tau)$ are  
invertible for $|\tau|\ge R$, $\tau \in \GR^*$, and some large $R$.  
In particular, by restricting the family $\hat{A}$ to the ball  
$$   
        B_R := \{ |\tau| \le R \},   
$$  
we obtain a family of elliptic operators that are invertible on the  
boundary of $B_R$, and hence $\hat{A}$ defines an element   
\begin{equation} \label{eq.def.degree}  
        \deg_{\GR}(A): = \deg(\hat{A}) \in K^0(B_R,\pa B_R) \simeq  
        K^0(\lgg^*)  
\end{equation}   
in the $K$-group of the ball of radius $R$, relative to its boundary,  
as explained above, called also the {\em degree} of $A$.  
  
We want a formula for the Chern character of the degree of   
$A \in \PsS {\infty}{Y;E_0,E_1}$.  
Because $A$ is elliptic, its principal symbol  
defines a class $[a'] \in K^0((T^*_{vert}Y)/\GR)$. If $E_0$ and  
$E_1$ are isomorphic, then it also defines a class   
$[a] \in K^0((S^*_{vert}Y)/\GR)$.   
Denote by $n$ the dimension of the quotient $Z_b=Y_b /\GR_b$ (which is  
independent of $b$ because we assumed $B$ connected), and  
let $\pi: (S^*_{vert}Y)/\GR \to B$ be the projection and  
\begin{equation}  
\begin{gathered}  
	\tilde\pi_* : H^*((T^*_{vert}Y)/\GR) \to H_c^{* - 2n }(\lgg^*)  
	\\ \pi_* : H^*((S^*_{vert}Y)/\GR) \to H_c^{* - 2n + 1}(\lgg^*)  
\end{gathered}  
\end{equation}  
be the integration along the fibers in cohomology.  Then  
$$  
	Ch(\deg_\GR(A)) \in H_c^*(\lgg^*)   
	\simeq H^{* +n}(B, \mathcal O),  
$$  
and the following theorem gives a formula for this cohomology class in  
terms of the classes $[a']$ or $[a]$ defined above.  
  
We denote by ${\mathcal T}$ the Todd class of the vector bundle  
$(T_{vert}Y)/\GR \otimes \CC \to Y/\GR$. We assume $B$ to be compact.

\begin{theorem}\label{Theorem.degree}\   
If $A \in \PsS{\infty}{Y;E_0,E_1}$ is elliptic, then the Chern  
character of $\deg_\GR(A)$ is given by  
\begin{equation*}  
	Ch( \deg_\GR(A)) = (-1)^n \tilde\pi_* \big (Ch[a'] {\mathcal T}\big)   
	\in H_c^{*}(\lgg^*),   
\end{equation*}  
Moreover, $Ch( \deg_\GR(A)) = (-1)^n \pi_* \big (Ch[a] {\mathcal  
T}\big)$ if $E_0 \cong E_1$.  
\end{theorem}

{\em Observations.} It is almost always the case that $E_0 \cong E_1$.  
For example, it is easy to see that this must happen if the Euler  
characteristic of $(T^*_{vert}Y)/\GR$ vanishes. Since $\GR = B\times  
\RR^q$ in many applications, this assumption is satisfied if $q > 0$.\  
  
Another observation is that if the set of elliptic elements in  
$\PsS{\infty}{Y;E_0,E_1}$ is not empty, $Z : = Y/\GR$ is compact.  
\vspace*{2mm}

\begin{proof}\ We cannot use Theorem \ref{Theorem.mult.b} directly  
because our family $\hat{A}$ does not consist of multiplication  
operators on the boundary. Nevertheless, we can deform $\hat{A}$ to a  
family of operators that are multiplication operators at $\infty$, for  
suitable $A$. We now construct this deformation.  
  
Let $E$ be a vector bundle over $V$. We consider classical symbols  
$S^m_c(E)$ whose support projects onto a compact subset of $V$.  Let  
$$  
	A_Y := (T_{vert}Y)/\GR \cong T_{vert}Z \times_B \lgg.  
$$  
First we need to define a nice quantization map $q : S^{m}(A_Y^*) \to  
\Psm{m}Y$. To this end, we proceed as usual, using local coordinates,  
local quantization maps, and partitions of unity, but being careful to  
keep into account the extra structure afforded by our settings:\ the  
fibration over $B$ and the action of $\GR$. Here are the details of  
how this is done.  
  
Fix a cross section for $Y \to Z := Y/\GR$ and, using it, identify $Y$  
with $Z \times_B \GR$ as $\GR$-spaces. Denote by $p_0 : Z \to B$ the  
natural projection. We cover $B$ with open sets $U_{\alpha'}$ that are  
diffeomorphic to open balls in a Euclidean space such that $Z  
\vert_{U_{\alpha'}} \cong U_{\alpha'} \times F$ and $\GR  
\vert_{U_{\alpha'}} \cong U_{\alpha'} \times \RR^q$. We also cover $F$  
with open domains of coordinate charts $V_{\alpha''}$.  Then we let  
$W_{\alpha} = p_0^{-1}(U_{\alpha'}) \cap V_{\alpha''}$, with $\alpha =  
(\alpha', \alpha'')$. The natural coordinate maps on $W_{\alpha}$ then  
give rise to a quantization map  
\begin{equation}  
	q_\alpha: S_c^m(A^*_Y\vert_{W_\alpha}) \to  
	 \Psm{m}{W_\alpha \times_B \GR}, \quad   
	q_\alpha(a) = a(b,x,D_x,D_t),  
\end{equation}  
where we identify  
\begin{equation}  
	S_c^m(A^*_Y\vert_{W_\alpha}) = S_c^m(T^*_{vert}(W_\alpha) \times_B \GR) =  
	S_c^m(U_{\alpha'} \times T^*V_{\alpha''} \times {\RR^*}^q).  
\end{equation}  
Denote by $b \in U_{\alpha'}$, $(x,y) \in T^*V_{\alpha''} \cong  
V_{\alpha''} \times {\RR^*}^n$, and $\tau \in {\RR^*}^q$ the  
corresponding coordinate maps. Then $q_\alpha(a) = a(b,x,D_x,D_t)$  
acts on $\CIc(U_{\alpha'} \times T^*V_{\alpha''} \times \RR^q)$ as  
\begin{multline*}  
	a(b,x,D_x,D_t)u(b,x,t) \\ = (2\pi i)^{-n - q} \int_{\RR^{n +  
	q}} \left (\int_{\RR^{n + q}}e^{i(x - z)\cdot y + i (t -  
	s)\cdot \tau} a(b,x,y,\tau) u(b,z,s) ds dz \right ) dy d\tau.  
\end{multline*}   
Choose now a partition of unity $\phi_\alpha^2$ subordinated to  
$W_\alpha$ and let  
$$  
	q(a) = \sum_\alpha q_\alpha(\phi_\alpha a) \phi_\alpha.  
$$   
  
The main properties of $q$ are the following:  
  
\begin{enumerate}  
\item\ if $a$ has order $m$, then $\sigma_m(q(a)) = a$, modulo symbols of  
lower order;  
\item\ there exist maps $q_b : S^m(T^* Z_b ) \to \Psi^m(Z_b)$ such  
that  
$$  
	\widehat{q(a)}(\tau) = q_b(a(\cdot,\tau)) = a(b,x,\tau + D_x),  
$$  
if $\tau \in \lgg_b^*$ and $x \in Z_b$;  
\item\ the maps $q_b$ define a quantization map   
$$  
	\tilde q : \Gamma( \lgg^*, S^m( T^*_{vert}Z)) \to  
	\Gamma(\lgg^*, \Psi^m(Z)),  
$$  
where we regard $\Psi^m(Z_b)$ as defining a bundle of algebras,  
$\Psi^m(Z)$, on $B$, first, and then on $\lgg^*$, by pull-back.  
Similarly, we regard $S^m(T^*_{vert}Z)$ as defining a bundle over $Z$,  
which we then pull back to a bundle on $\lgg^*$.  (See also the  
discussion related to Equations \eqref{eq.bundle.hatA} and  
\eqref{eq.section.hatA}.)  
\end{enumerate}  
  
The deformation of our family is obtained as follows. Let $|\tau|$ and 
$|y|$ be the norms $\tau \in \lgg^*$ and $y \in T^*_{vert}Z$. Define 
then 
\begin{equation}  
	\phi_\lambda^2 (y,\tau) = 1 + \lambda |y|^2 + \lambda (1 + 
	\lambda |\tau|^2)^{-1} |y|^2 \quad \text{and }\; 
	\psi_\lambda^{-2} (y,\tau) = 1 + \lambda |\tau|^2, 
\end{equation}  
which are chosen to satisfy $1 + \phi_\lambda^2 |\tau|^2 +  
\psi_\lambda^2 |y|^2 = 1 + |\tau|^2 + |y|^2 + \lambda |\tau|^2|y|^2$.  
For any symbol $a \in S^0(A_Y^*)$, $A_Y^* = T^*_{vert}Z \times_B  
\lgg^*$, we let  
$$  
	a_{\lambda,\tau}(y) = a( \psi_\lambda y, \phi_\lambda \tau ),  
	\quad \lambda \in [0,1],\; \tau \in \lgg^*_b, \text{ and } y  
	\in T^*Z_b.  
$$  
We can define then $A_{\lambda}(\tau) := q_b(a_{\lambda,\tau})$, $\tau  
\in \lgg^*_b$, which is the same as saying that $A_\lambda =  
\tilde{q}(a_{\lambda,\tau})$, and these operators will satisfy the  
following properties:  
  
\begin{enumerate}   
\item\ For each fixed $\lambda$, the operators $A_{\lambda}(\tau)$  
define a section of $\Psi^m(Z)$ over $\lgg^*$ and these sections  
depend smoothly on $\lambda$ (in other words, $A_{\lambda}(\tau)$  
depends smoothly on both $\lambda$ and $\tau$, in any trivialization);  
\item\ $A_{0}(\tau) = \widehat{q(a)}(\tau)$, for all $\tau$;  
\item\ For each nonzero $\tau' \in \lgg^*$ and $\lambda >0$, the limit  
$\displaystyle{\lim_{t \to \infty}} A_{\lambda}(t\tau')$ exists and is  
a multiplication operator;  
\item If $a,b \in S^0(A_Y^*)$ are such that $ab = 1$, $a$ is  
homogeneous of order zero outside the unit ball, and if we define  
$A_{\lambda} := \tilde{q}(a_{\lambda,\tau})$ and $B_{\lambda,\tau} :=  
\tilde{q}(b_{\lambda,\tau})$, then there exists a constant $C >0$ such  
that  
$$  
	\| A_\lambda(\tau) B_\lambda(\tau) - 1 \| \le C(1 + |\tau|)^{-1}  
$$  
and similarly  
$$  
	\| B_\lambda(\tau) A_\lambda(\tau) - 1 \| \le C(1 + |\tau|)^{-1}  
$$  
For all $\tau$ and $\lambda$;  
\item All these estimates extend in an obvious way to matrix valued  
symbols.  
\end{enumerate}  
  
These properties are proved as follows. We first recall that, for any 
vector bundle $E$, we can identify the space of classical symbols 
$S_c^0(E)$ with $\CIc(E_1)$, the space of compactly supported 
functions on $E_1$, the unit ball of $E$, by $E_1 \smallsetminus \pa 
E_1 \ni \xi \to (1 - \|\xi\|^2)^{-1}\xi \in E$. For any quantization 
map, the norm of the resulting operator will depend on finitely many 
derivatives. Because we can extend $a_{\lambda,\tau}$ to a smooth 
function on the radial compactification of $A_Y^*$, the first property 
follows. The second property is obvious.  The third property is 
obtained using the same argument and observing that, for $\lambda >0$, 
we can further extend our function to the radial compactification in 
$\lambda$ also. By investigating what this limit is along various 
rays, we obtain the third property. 
  
The fourth property is obtained using the following  
observation: there exist a constant $C>0$ and seminorms $\|\,\cdot\,\|_0 $ on  
$S_c^0(T^*\RR^n)$ and $\|\,\cdot\,\|_{-1}$ on $S^{-1}(T^*\RR^n)$  
such that, for any symbols $a,b \in S_c^0(T^*\RR^n)$,  
$$  
	\|(ab)(x,D_x) - a(x,D_x) b(x,D_x) \| \le \sum_j    
	C(\|a\|_0 \| \pa_{y_j} b\|_{-1}  
	+ \| \pa_{y_j} a\|_{-1}\|b\|_0 ),  
$$  
$\pa_{y_j}$ being all derivatives in the symbolic directions (whose  
coordinates are denoted by $y$).   
Finally, the fifth property is obvious.  
  
We now turn to the proof of the formula for the degree of $A$ stated  
in our theorem. We prove it by a sequence of successive reductions,  
using the facts established above. First, it is easy to see that  
$\deg_\GR(A)$ depends only on its principal symbol, and hence we can  
assume that $A$ has order zero and $A = q(a)$, where $a =  
\sigma_0(A)$.  
  
The above deformation can be used to prove our theorem as follows.  
Fix $R$ large enough, and restrict the families $A_{\lambda}$ to the  
closed ball of radius $R$ in $\lgg^*$. For $|\tau| = R$ large enough,  
all operators $A_{\lambda}(\tau)$, $\lambda \in [0,1]$ are invertible,  
so the degree $\deg (A_\lambda)$ of these families is defined and does  
not depend on $\lambda$ or $R$, provided that $R$ is large  
enough. Since $\deg_\GR(A) = \deg(A_0)$, by definition, it is enough  
to compute $\deg(A_\lambda)$, for any given $\lambda$.  Choose then  
$\lambda >0$ arbitrary, and let $R \to \infty$.  Then the family  
$A_\lambda$ extends to a continuous family on the radial  
compactification of $\lgg^*$, which consists of multiplication  
operators on the boundary. Moreover, the symbol class of $A_\lambda$  
is nothing by the extension of $a_{\lambda}(\tau)$ to the radial  
compactification in $\tau$ and $\lambda$ (which is a manifold with  
corners of codimension two).  
  
We can use then Theorem \ref{Theorem.mult.b} to conclude that   
$$  
	\deg(A_\lambda) = (-1)^{n} \tilde\pi_* ( Ch [a_\lambda']  
	{\mathcal T}) \in H^*_c(\lgg^*).  
$$  
But $a_\lambda$ is homotopic to $a$ through symbols that are  
invertible outside a fixed compact set, so $[a_\lambda']=[a']$.  We  
get  
$$  
	\deg(A) = (-1)^{n} \tilde\pi_* ( Ch [a']  
	{\mathcal T}) \in H^*_c(\lgg^*).  
$$  
  
To obtain the second form of our formula for $E_0 \cong E_1$,  
and thus finish the proof, we proceed as at the end of the proof of  
Theorem \ref{Theorem.mult.b}, using $Ch [a'] = \pa Ch [a]$.    
\end{proof}

To prove the following result, we shall use terminology from algebraic  
topology:\ if $I_k \subset A_k$ are two-sided ideal of some algebras  
$A_0$ and $A_1$ and $\phi:A_0\to A_1$ is an algebra morphism, we say  
that $\phi$ {\em induces a morphism of pairs} $\phi:(A_0,I_0) \to  
(A_1,I_1)$ if, by definition, $\phi(I_0) \subset I_1$.

\begin{theorem}\label{Theorem.inddeg}\ Let $\GR$ be a bundle of   
abelian Lie groups and $A \in \PsS m{Y;\CC^N}$ be an elliptic  
element. Then  
$$   
        \ind_a(A) = \deg_{\GR}(A) \in K^0(\lgg^*).   
$$   
\end{theorem}

\begin{proof}\  
Let $B_R = \{ |\tau| \le R\} \subset \GR^*$ be as above.  The  
algebra   
$$  
      {\mathfrak A}_R := \CI(B_R,\Ps {\infty}{Y_b})  
$$  
of $\CI$-families of pseudodifferential operators on $B_R$ acting on  
fibers of $Y \times_B B_R \to B_R$, contains as an ideal ${\mathfrak  
I}_R = \CI_0(B_R,\Ps {-\infty}{Y_b})$, the space of families of  
smoothing operators that vanish to infinite order at the boundary of  
$B_R$.  If $A$ is an elliptic family, as in the statement of the  
lemma, and if $R$ is large enough, then $\hat{A}$, the indicial family  
of $A$, defines by restriction an element of $M_N({\mathfrak A}_R)$  
that is invertible modulo $M_N({\mathfrak I}_R)$.  
  
Recall that the boundary map $\pa_1$ in algebraic $K$-theory  
associated to the ideal ${\mfk I}_R$ of the algebra ${\mathfrak A}_R$  
gives $ \pa_1[A]= \deg_{\GR}(A),$ by definition. Also, the boundary  
map $\pa_0$ in algebraic $K$-theory associated to the ideal  
$\PsS{-\infty}Y$ of the algebra $\PsS {\infty}{Y}$ gives  
$\pa_0[A]=\ind_a(A)$. We want to prove that $\pa_1[A]=\pa_0[A]$.  The  
desired equality will follow by a deformation argument, which involves  
constructing an algebra smoothly connecting the ideals ${\mfk I}_R$  
and $\PsS {-\infty}Y$.  
  
Consider inside $\CI([0,R^{-1}], \PsS{-\infty}Y)$ the subalgebra of  
families $T=(T_x)$ such that $\hat{T}_x(\tau) = 0$ for $|\tau| \ge  
x^{-1}$. (In other words, $T_x \in {\mfk I}_{x^{-1}}$, if $x \not =  
0$, and $T_{0}$ is arbitrary.)  Denote this subalgebra by ${\mfk  
I}_{R\infty}$. Also, let ${\mfk A}_{R\infty}$ be the set of families  
$A=(A_x)$, $x \in [0,R^{-1}]$, $A_x \in {\mfk A}_{x^{-1}}$, if $x \not  
=0$, $A_0 \in \PsS {\infty}Y$ arbitrary such that the families $AT :=  
(A_xT_x)$ and $TA =: (T_x A_x)$ are in ${\mfk I}_{R\infty}$, for all  
families $T=(T_x) \in {\mfk I}_{R\infty}$.  
  
It follows that ${\mfk I}_{R\infty}$ is a two-sided ideal in ${\mfk  
A}_{R\infty}$ and that the natural restrictions of operators to  
$x=R^{-1}$ and, respectively, to $x=0$, give rise to morphisms of  
pairs  
$$  
       e_{1}:({\mfk A}_{R\infty}, {\mfk I}_{R\infty}) \to   
       ({\mfk A}_{R}, {\mfk I}_{R}) , \quad \text{ and }  
$$  
$$  
       e_{0}: ({\mfk A}_{R\infty}, {\mfk I}_{R\infty}) \to (\PsS  
       {\infty}Y, \PsS {-\infty}Y).  
$$  
Moreover, the indicial family of the operator $A$ gives rise, by  
restriction to larger and larger balls $B_r$, to an invertible element  
in ${\mfk A}_{R\infty}$, also denoted by $A$. Let $\pa$ be the  
boundary map in algebraic $K$-theory associated to the pair $({\mfk  
A}_{R\infty},{\mfk I}_{R\infty})$. Then $(e_0)_* \pa[A]= \pa_0[A]$ and  
$(e_1)_* \pa[A]=\pa_1[A]$. Since $(e_0)_*: K_0({\mfk I}_{R\infty}) \to  
K_0(\PsS {-\infty}Y)$ and $(e_1)_*: K_0({\mfk I}_{R\infty}) \to  
K_0({\mfk I}_{R})$ are natural isomorphisms, our result follows.  
\end{proof}

We now drop the assumption above that $\GR$ consist of abelian Lie 
groups, assuming instead that $\GR$ consists of simply-connected 
solvable Lie groups, and want to compute the Chern character of the 
analytic index $\ind_a(A)$, for an elliptic family $A \in \Psm {m}Y$. 
One difficulty that we encounter is that the space on which the 
principal symbols are defined, that is $(S^*_{vert}Y)/\GR$, is not 
orientable in general. (Recall that $S^*_{vert}Y$ is the space of 
vectors of length one of $T_{vert}^*Y$, the dual of the vertical 
tangent bundle $T_{vert}Y$ to the fibers of $Y \to B$.) 
  
We denote by ${\mathcal T}$ the Todd class of the vector bundle  
$(T_{vert}Y)/\GR \otimes \CC \to Y/\GR$ and by $\pi_*$ the integration  
along the fibers of $(S^*_{vert}Y)/\GR \to B$, as above. We assume $B$   
to be compact.

\begin{theorem}\label{Theorem.Chern}\   
Let $\GR$ be a bundle of Lie groups whose fibers are simply-connected,  
solvable Lie groups. Let $A \in \Psm {m}{{Y,\CC^N}}$ be an elliptic,  
invariant family, and let $[\sigma_m(A)] \in K^1((S^*_{vert}Y)/\GR)$ be  
the class defined by the principal symbol $\sigma_m(A)$ of $A$. Then  
the Chern character of the analytic index of $A$ is given by  
$$   
        Ch(\ind_a(A)) = (-1)^n \pi_*\big (Ch[\sigma_m(A)]   
        {\mathcal T}\big) \in H_c^{*}(\lgg^*),  
$$  
where $n$ is the dimension of the fibers of $(S^*_{vert}Y)/\GR \to B$.  
\end{theorem}

\begin{proof}\   
Note first that we can deform the bundle of Lie groups $\GR$ to the  
bundle of {\em commutative} Lie groups $\lgg$ as before, using  
$\GR_{ad}$. Moreover, we can keep the principal symbol of $A$ constant  
along this deformation. This shows that we may assume $\GR$ to consist  
of commutative Lie groups, i.e. that $\GR$ is a vector  
bundle. The result then follows from Theorems \ref{Theorem.degree} and  
\ref{Theorem.inddeg}.  
\end{proof}

{\em Observations.} We can extend the above theorems in several ways.  
First, we can drop the assumption that $Z= Y/\GR$ be compact, but then  
we need to consider bounded, elliptic elements $A \in \Hom(E_0,E_1) +  
\Psm 0{Y;E_0,E_1}$ (or, if $\GR$ consists of abelian Lie groups, then 
$A \in \Hom(E_0,E_1) + \Psm 0{Y;E_0,E_1}$).  Also, in the last two theorems, we 
can allow operators acting between sections of different vector 
bundles. This will require to slightly modify the proof of Theorem 
\ref{Theorem.inddeg}, either by using a smooth version of bivariant 
$K$-theory \cite{NistorKthry}, or by using the usual bivariant 
$K$-theory after we have taken the norm closures of the various ideals 
$\mfk I$ decorated with various indices. We can also further integrate 
along the fibers of $\lgg^* \to B$ to obtain a cohomological formula 
with values in $H_c^*(B;\mathcal O)$, the cohomology with local 
coefficients in the orientation sheaf of $\lgg^* 
\to B$. This will be useful in Section~\ref{S.HI}.

\section{Regularized traces\label{S.RT}}

Having in mind future applications, we also want to give a local  
formula for the equivariant family index of an invariant, elliptic  
family of operators, as considered in the previous section.  This will  
be done in terms of various residue type traces.  In this section, we  
develop the analytic tools required to define these regularized  
traces.\\[2mm]  
{\bf Assumption.}\ {\em Throughout the rest of this paper, we shall 
assume that $\GR$ consist of simply-connected, non-trivial abelian Lie 
groups, and hence that it is a non-zero vector bundle.} 
  
\vspace*{2mm}Recall that $Y$ is a flat $\GR$-bundle if $Y = B \times Z  
\times \RR^q$ and $\GR = B \times \RR^q$ ($q > 0)$. The results we  
will establish are local in $B$, and hence we can reduce the general  
case to the flat case. Actually, it is easier to assume first that $B$  
is reduced to a point. We thus carry the analysis first in this case,  
and then we extend the results to the general case. The Lemma  
\ref{lemma.inv.t} and Proposition \ref{prop.ext} are probably not new.  
We nevertheless include their proofs for completeness and to fix notation.  
  
There is an action of $\RR^q$ on $\PsS {-\infty}{Z \times \RR^q}$, the  
action of $\xi \in \RR^q$ is obtained by multiplying the convolution  
kernel of an operator $A \in \PsS p{Z \times \RR^q}$ by  
$\exp(it\cdot\xi),$ where $t \in \RR^q$ are coordinates for the second  
component in $Z \times \RR^q$. In terms of the Fourier transform  
representation of these operators, the action of $\xi$ becomes  
translation by $\xi\in\RR^q.$  
  
We shall denote by $Tr$ the usual (Fredholm) trace on the space of  
trace class operators on a given Hilbert space.

\begin{lemma}\label{lemma.inv.t}\   
The space of $\RR^q$--invariant traces on $\PsS {-\infty}{Z \times  
\RR^q}$ is one-dimensional.  
\end{lemma}

\begin{proof}\ Consider the map  
\begin{equation}  
        \aTr(A)=\int Tr(\hat A(\tau))\,d\tau.  
\label{eq.13}\end{equation}  
We need to show that this is the only invariant trace functional.  In  
terms of indicial families, the infinitesimal generators of the  
$\RR^q$--action correspond to the multiplication operators with the  
functions $t_k$.  Let  
$$  
        \Hd_0(\PsS {-\infty}{Z \times \RR^q}) :=   
        \PsS {-\infty}{Z \times \RR^q} / [\PsS {-\infty}{Z \times \RR^q},  
        \PsS {-\infty}{Z \times \RR^q} ]  
$$  
be the first homology group of $\PsS {-\infty}{Z \times \RR^q}$.  It  
remains to show that the subspace of $\Hd_0(\PsS {-\infty}{Z \times  
\RR^q})\simeq {\mathcal S}(\RR^q)$ spanned by $t_k\Hd_0(\PsS  
{-\infty}{Z \times \RR^q})$ has codimension $1$.  Indeed, the kernel  
$\Hd_0(\PsS {-\infty}{Z \times \RR^q}) \to \CC$ of the evaluation at  
$0$ is the span of $t_k\Hd_0(\PsS {-\infty}{Z \times \RR^q})$.  This  
proves the lemma.  
\end{proof}

Let $x$ be the identity function on $[0,\infty)$ and $l_s$  
be a smooth function on $[0,\infty)$ such that  
$l'_s(x)=x^{s-1}$ for $x \ge 1$. (So that, in particular, $l_{0}(x) = \ln x$,  
for $x$ large.) We define then the spaces of functions   
$$  
        {\mathcal M}_k=S^{\infty}([0,\infty)) + \CC[x]l_0, \quad \mbox{ for }\,  
        k \in \ZZ,  
$$  
and  
$$  
        {\mathcal M}_s = S^{\infty}([0,\infty))l_s, \quad \mbox{ for }  
        s \in \CC \smallsetminus \ZZ.  
$$   
Thus, $\mathcal{M}_0$ consists of smooth functions on $[0,\infty)$,  
that can be written, for any $M \in \ZZ_+$, as  
\begin{equation}\label{eq.zero.space}  
        f(x)=h_M(x) + \sum_{k = -M}^{-1}a_k x^k  + \sum_{k=0}^N   
        (b_k + c_k \log x)x^k, \quad \forall x \ge 1,   
\end{equation}  
where $a_k,b_k,c_k$ are complex parameters, $N\in \ZZ_+$, and  
$h_M \in S^{-M-1}([0,\infty))$. Similarly, the  
space $\mathcal{M}_s$, $s\not \in \ZZ$,  
consists of smooth functions $f\in \CI([0,\infty))$ that  
can be written, for any $M \in \ZZ_+$, as  
\begin{equation}\label{eq.s.space}  
        f(x)=h_M(x)x^s + \sum_{k=-M}^{N}\alpha_k x^{k+s},  
        \quad\forall x \ge 1.  
\end{equation}  
for some constants $\alpha_k \in \CC$, $N \in \ZZ$, and $h_M \in  
S^{-M-1}([0,\infty))$. Fix $M > \max\{s,0\}$ and $R \ge 1$, and define  
\begin{gather}  
        I(f)(x)=\int_0^x f(x)dx\ , \\  
        \ITr f =\int_0^R  f(x)dx + \int_R^\infty h_M(x)dx  
        - \sum_{k=-M}^{-2}\frac{a_kR^{k+1}}{k+1} - a_{-1} \log R\\  
        - \sum_{k=0}^{N}\frac{R^{k+1}}{k+1}  
        \big(b_k - \frac{c_k}{(k+1)} +c_k \log R \big), \quad   
        \Mif f \in {\mathcal M}_l, \  l \in \ZZ  
\end{gather}  
and  
\begin{multline}\label{eq.(24)} 
        \ITr f(x)dx = \ITr f=\int_0^R  f(x)dx + \int_R^\infty h_M(x)x^sdx  
        - \sum_{k=-M}^{N}\frac{\alpha_kR^{k+s+1}}{k+s+1}, \\  
        \Mif f \in {\mathcal M}_s,\   
        s \not \in \ZZ.  
\end{multline}  
It is easy to see that $I$ maps $\mathcal{M}_0$ to itself   
and that the definition of $\ITr f$ is independent of $M$ and $R$,  
for $f \in \mathcal{M}_s$, $s \in\CC$. Moreover,   
\begin{equation}\label{eq.inv.derivative}  
        \ITr I^k(f)=\ITr I^{k+1}(f'),  
\end{equation}  
if $f \in \mathcal{M}_0$ and $k \ge 0$, because $I^k(f)-I^{k+1}(f')$  
is a polynomial of degree at most $k$ and $\ITr$ vanishes on  
polynomials. If $s\not \in \ZZ$, then $I$ maps $\mathcal{M}_s$ to  
$\mathcal{M}_s + \CC[x]\simeq \mathcal{M}_s \oplus \CC[x]$, and since  
this is a direct sum decomposition, we can extend $\ITr$ to  
$\mathcal{M}_s + \CC[x]$ to vanish on polynomials, which guaranties  
that the Equation \eqref{eq.inv.derivative} is still  
satisfied. Moreover, if $f \in \mathcal{M}_s$ is as in Equation  
\eqref{eq.s.space}, $P(x) = \sum_{k=0}^N b_kx^k$, and $g = f + P$,  
then the original definition of $\ITr$ is still valid in this case,  
with the obvious changes:  
\begin{equation}  
        \ITr g=\int_0^R g(x)dx + \int_R^\infty h_M(x)x^sdx -  
	\sum_{k=-M}^{N}\frac{\alpha_kR^{k+s+1}}{k+s+1} -  
	\sum_{k=0}^{N}\frac{b_kR^{k+1}}{k+1}.  
\end{equation}  
  
The following lemma shows that $\ITr$ provides a natural extension of
the integral on $[0,\infty)$ (which justifies the notation $\ITr f =
\ITr_0^\infty f(x) dx$), in the following sense.

\begin{lemma}\label{lemma.hol.reg}\ 
Let $f_z(x) \in {\mathcal S}^{m} ([0,\infty))$ be holomorphic in $z$
and $g_z(x) := f_z(x) (1 + x)^{-z}$. Then $\ITr_0^\infty g_z(x)dx$ is
holomorphic on $\CC \smallsetminus (m + \ZZ)$, with at most simple
poles at $m + \ZZ$, and $\ITr_0^\infty g_z(x)dx = \int_0^\infty g_z(x)
dx$ for $\Re(z) > m +1$, where the second integral is absolutely
convergent.
\end{lemma}

\begin{proof}\ This follows from the fact that the right
hand side of Equation \eqref{eq.(24)} is holomorphic in $f$ and $s$
for $s \not \in \ZZ$. The same formula guarantees at most simple poles
at $s \in \ZZ$. Since $s = z - m$, in our case, the result follows.
\end{proof}

If $f \in \CI(\RR)$ is such that $f_+, f_- \in {\mathcal M}_s +  
\CC[X]$, for some $s$, where $f_+(\tau) = f(\tau)$ and  
$f_-(\tau) = f(-\tau)$, $\tau \ge 0$, we define  
\begin{equation*}  
  	\ITr f = \ITr f_+ + \ITr f_-.
\end{equation*}  
  
{\em We fix from now on a positive, invertible operator $D \in \PsS
1{Z \times
\RR^q}$.}  Let $\CC \ni z \to A(z) \in \PsS {m}{Z \times \RR^q}$ be an 
entire function.  Then 
\begin{equation}\label{eq.f.z}  
        f_z(\tau) = Tr \big [\pa_\tau^\alpha \big(|\tau|^k
        \hat{D}(\tau)^{-z} \hat{A} (z,\tau)\big) \big ], \quad k \in
        \ZZ_+
\end{equation}   
is defined and holomorphic for any multi-index $\alpha$, for any
$\Re(z -m)> \dim Z - |\alpha| = d - |\alpha|$, and for any fixed $\tau
\in \RR^q$.  Moreover, by classical results (see
\cite{Gilkey1}, for example), the function $z \to f_z(\tau)$ has a 
meromorphic extension to $\CC$, for each fixed $\tau$, with at most 
simple poles at integers. Let  
\begin{equation*} 
	\Omega = (\CC \smallsetminus \ZZ) \cup \{z \in \CC ,\Re(z - m) > 
	d - |\alpha| \}. 
\end{equation*} 
(Recall that $d = \dim Z$). In the following proposition $m$ does not
have to be an integer.

\begin{proposition}\label{prop.ext}\   
Let $A(z) \in \PsS {m}{Z \times \RR^q}$ be an entire function and $D
\in \PsS {1}{Z \times \RR^q}$ be an invertible, positive operator.
Also, let $f_z(\tau)$, defined for $z \in \Omega$, be as in Equation
\eqref{eq.f.z} above. 
  
(i) The function $f_z(\tau)$ is in $\CI(\Omega \times \RR^q)$ and the  
map $z \to f_z(\tau)$ is holomorphic on $\Omega$, for each fixed $\tau  
\in \RR^q$.  
  
(ii) There is a holomorphic $g:\Omega \to S^{m+d-|\alpha|}(\RR^q)$ such that 
$f_z(\tau)=g_z(\tau)|\tau|^{k-z}$, for all $\tau$ such that $|\tau| 
\ge 1$, and hence $f_z(x\tau) \in \mathcal{M}_{m+d+k-z}$ (as a 
function of $x$), for each fixed $\tau \not = 0$. 
  
(iii) The function $z \to \ITr_0^\infty f_z(x\tau)dx$ is holomorphic
on $\CC \smallsetminus (m + \ZZ)$, with at most simple poles at $m +
\ZZ$.
\end{proposition}

\begin{proof}\     
The proof for $k \not = 0$ or $\alpha \not = (0,\ldots,0)$ is the same
as that for $k = 0$ and $\alpha = (0,\ldots,0)$, so we shall assume
that we are in the latter situation. Also, by replacing $z$ by $z-m$,
we can assume that $m \in \ZZ$.
  
We first prove the lemma for $m = -\infty$, that is for $A(z) \in \PsS  
{-\infty}{Z \times \RR^q}$.  Denote by $\mathcal{K}$ the algebra of  
compact operators acting on $L^2(Z)$ and by $\mathcal{C}_1  
\subset\mathcal{K}$ the normed ideal of trace class operators. For any  
$M \in \ZZ_+$, the product $\hat{D}^M(\tau)\hat{A}(z,\tau)$ is in  
$\mathcal{S}(\RR,\mathcal{C}_1) =\mathcal{S}(\RR) \widehat{\otimes}  
\mathcal{C}_1$ (here $\widehat{\otimes}$ denotes the completed  
projective tensor product).  Also, because $D$ is invertible and  
positive, the function $(z,\tau) \to \hat{D}(\tau)^{-z} \in  
\mathcal{K}$ is differentiable, with bounded derivatives, and  
holomorphic in $z$, for $\Re(z) \ge 1$. Since $\Tr: \mathcal{K}  
\hotimes \mathcal{C}_1 \to \CC$ is continuous, it follows that the  
function  
\begin{equation*}  
        (z,\tau) \to Tr(\hat{D}(\tau) ^{-z-M-1} 
        \hat{D}(\tau)^{M+1}\hat{A}(z,\tau)) \in \CC 
\end{equation*}  
is differentiable, with bounded derivatives, and holomorphic in $z$  
for $\Re(z) \ge -M$. Since $M$ is arbitrary, this proves (i) and (ii)  
for $A(z) \in \PsS {-\infty}{Z \times \RR^q}$. The last statement is  
an immediate consequence of (ii), because $z \to g_z \in  
\mathcal{S}(\RR^q)$ is holomorphic.  
  
Using now the fact that the lemma is true for $A(z)$ in the residual  
ideal, we may assume, using a partition of unity, that $Z=\RR^d$ and  
that the Schwartz convolution kernels of $\hat{A}(z,\tau)$ are  
contained in a fixed compact set.  
  
Let $\Delta_0,\Delta_1\ge 0$ be the constant coefficient Laplacians on   
$Z=\RR^d$ and $\RR^q$, respectively. We define  
$D_0=(1 + \Delta_0 + \Delta_1)^{1/2} \in \PsS {1}{Z \times \RR^q}$.   
To prove the lemma for $A(z) \in \PsS {m}{Z \times \RR^q}$, $m >\infty$,   
we shall first assume that $D = D_0$. Clearly, $D_0 \in   
\PsS {1}{Z \times \RR^q}$. If $A(z)=a(z,x,D_x,D_t)$, for a symbol   
$a(z,\cdot,\cdot,\cdot) \in S^{m}(T^* Z \times \RR^q)=S^m(\RR^d \times \RR^d  
\times \RR^q)$, then   
$$  
        \hat{A}(z,\tau) = a(z,x,D_x,\tau) \quad \text{ and } \quad   
        \hat{D}_0(\tau)^{-z}=(1 + \Delta_0 + |\tau|^2)^{-z/2}.  
$$   
This gives, by the standard calculus, that   
$\hat{A}(z,\tau)\hat{D}_0(\tau)^{-z}=a_1(z,x,D_x,\tau)$ for  
\begin{equation}  
        a_1(z,x,\xi,\tau)=a(z,x,\xi,\tau)(1+|\tau|^2+|\xi|^2)^{-z/2}.  
\end{equation}  
  
{}From the above relation, we obtain   
\begin{multline*}  
        f_z(\tau):=Tr(\hat D(\tau)^{-z}\hat A(\tau))
        =Tr\big(\hat{A}(\tau)\hat{D}(\tau)^{-z}\big)\\
        =(2\pi)^{-d}\int_{T^*Z}a(z,x,\xi,\tau)
        (1+|\tau|^2+|\xi|^2)^{-z/2} d\xi dx, \quad \Mfor \ \tau\ge 1 .
\end{multline*}  
Using the asymptotic expansion of $a$ in homogeneous functions in
$(\xi,\tau)$ and the substitution $\xi\to (1 + |\tau|^2)^{1/2}\xi$,
and the asymptotic expansion of $(1 + |\tau|^2)^{1/2}$ in powers of
$|\tau|$ at $\infty$, we obtain (i) and (ii) for this particular
choice of $D=D_0$.

We obtain (iii) directly from (ii) using Lemma \ref{lemma.hol.reg}.

The case $D$ arbitrary follows by writing
$D^{-z}=D_0^{-z}\big(D_0^{z}D^{-z} \big)$ and observing that $\CC \ni
z \to D_0^{z}D^{-z}\in\PsS {0}{Z \times \RR^q}$ is an entire function.
\end{proof}

See also \cite{KontsevichVishik}.  
 
We shall also need the following consequence of Proposition  
\ref{prop.ext}, above.

\begin{corollary}\label{cor.ext.one}\ 
Using the notation of the above proposition, we have that the  
function  
\begin{equation*} 
	F(s) := \int_{S^{q-1}} \left ( \ITr_0^\infty I^l \circ \pa_x^l
	\big[x^{q-1}f_z(x \tau) \big] dx \right ) d\tau
\end{equation*} 
is holomorphic on $ \{z \in \CC, \Re(z - m) > q + d\} \cup \CC
\smallsetminus (m + \ZZ)$, with at most simple poles at $m + \ZZ$, and
extends the function $\int_{\RR^q} f_z(\tau) d\tau$, which is defined
for $\Re(z - m) > q + d$.
\end{corollary}

\begin{proof}\ 
Assume $l=0$. The proof for arbitrary $l$ is completely similar.  The
function $\ITr_0^\infty x^{q-1}f_z(x \tau) dx$ is a holomorphic
extension of the function $\int_0^\infty x^{q-1}f_z(x \tau) dx$, which
is convergent for $\Re(z)$ large.  The result is obtained then by
integration in polar coordinates and by combining Lemma
\ref{lemma.hol.reg} with (ii) of the above proposition.
\end{proof}

Assume for the moment that $q = 1$ (and hence that $Y = Z \times \RR$).   
Using the above lemma and the functionals $\ITr$ and $I$, we obtain, as in 
\cite{Melrose46}, a functional $\bTr 1$ on $\PsS {\infty}{Z \times  
\RR}[|\tau|]$, by the formula  
\begin{equation}  
        \bTr 1(|\tau|^jA)=\ITr I^k(f_+ + f_-), 
\end{equation}  
where $f_+(\tau)= Tr \big [ \pa_{\tau}^k ( |\tau|^j\hat{A}(\tau)) \big
]$, $f_{-}(\tau)= Tr \big [ \pa_{\tau}^k ( |\tau|^j\hat{A}(- \tau))
\big ]$, if $\tau \le 0$, and $k \in \ZZ_+$, $k > m + \dim Z +
1$. From equation \eqref{eq.inv.derivative}, we see that this
definition is independent on $k$. The tracial property of $\bTr 1$
follows, as in \cite{Melrose46}, from
$$  
        \pa_\tau[\hat{A}(\tau),\hat{B}(\tau)]  
        =[\pa_\tau\hat{A}(\tau),\hat{B}(\tau)]   
        + [\hat{A}(\tau),\pa_\tau\hat{B}(\tau)].  
$$  
  
Let now $q$ be arbitrary, but we continue to assume that $B$ is reduced  
to a point. The following lemma will allow us to generalize the definition 
of $\bTr 1$.

\begin{lemma}\label{Lemma.restr.iso}\  
Restriction of the indicial family $\hat{A}$ to $\RR x$, $x \in  
\SS^{q-1}$, defines an $O(q)$-equivariant family of algebra morphism  
$r_{x}:\PsS {\infty}{Z \times \RR^q} \to \PsS {\infty}{Z \times 
\RR}$. Each $r_x$ restricts to a degree preserving isomorphism $\PsS  
{\infty}{Z \times \RR^q}^{O(q)} \simeq \PsS {\infty}{Z \times  
\RR}^{\ZZ_2}$, which is independent of $x$.  
\end{lemma}

\begin{proof}\ It is clear that the restriction of $\hat A$ to a line  
$\RR x$ is in $\mathcal{S}(\RR x,\Psi^{-\infty}(Z))$, whenever $A$  
is in $\PsS {-\infty}{Z \times \RR^q}$.  Moreover, we obtain  
isomorphisms  
\begin{equation*}  
        \PsS {-\infty}{Z \times \RR^q}^{O(q)}   
        \simeq \mathcal{S}(\RR^q,\Psi^{-\infty}(Z))^{O(q)}  
        \simeq \mathcal{S}(\RR,\Psi^{-\infty}(Z))^{\ZZ_2}  
        \simeq \PsS {-\infty}{Z \times \RR}^{\ZZ_2}.  
\end{equation*}  
These isomorphisms allow us to assume, using a partition of unity  
argument,  that $Z=\RR^l$. Using the fact that a symbol  
$a \in S^m(T^*Z \times \RR^q)$ restricts to a symbol in   
$S^m(T^*Z \times \RR x)$, $x \in \SS^{p-1}$, and the relation   
\begin{equation*}  
        \hat A(\tau)=a(x,D_x,\tau),  
\end{equation*}  
if $A=a(x,D_x,D_\tau)$, we see that  
the restriction of $\hat{A}$ to $\RR x$ is the indicial family  
of an operator in $\PsS {m}{Z \times \RR}$, denoted $r_{x}(A)$. Since    
\begin{equation*}  
        S^m(T^*Z \times \RR^q)^{O(q)} \simeq S^m(T^*Z \times
        \RR)^{\ZZ_2},
\end{equation*} 
the isomorphism  $\PsS {\infty}{Z \times \RR^q}^{O(q)} \simeq    
\PsS {\infty}{Z \times \RR}^{\ZZ_2}$ follows.  
\end{proof}

Let $A \in \PsS {\infty}{Z \times \RR^q}$ and denote by  
$A_1=\int_{O(q)}v(A)dv$ its average over $O(q)$ with respect to its  
normalized Haar measure, which we identify with an element of $\PsS  
{\infty}{Z \times \RR}^{\ZZ_2}$, thanks to Lemma  
\ref{Lemma.restr.iso}. Define  
\begin{equation}\label{eq.Tr.q}  
        \bTr q(A)= \operatorname{vol}(S^{q-1}) \bTr  
        1(|\tau|^{q-1}A_1).  
\end{equation}

\begin{lemma}\label{Lemma.on.trace}\   
The functional $\bTr q$ is an $O(q)$--invariant trace on $\PsS  
{\infty}{Z \times \RR^q}$, which extends the trace $\aTr$ defined on  
$\PsS {-s}{Z \times \RR^q}$, $\Re(s) > d + q$ ($d = \dim Z$),  
by Equation \eqref{eq.13}.  
\end{lemma}

\begin{proof}\  
The map $\bTr q$ is obviously well defined and $O(q)$-invariant, in
view of the above lemma.  In order to check the tracial property, we
use the definition. Fix $x \in\RR^q$ of length one, arbitrarily, then
$$ \bTr q(A)= \operatorname{vol}(S^{q-1})\int_{O(q)}\bTr
1(|\tau|^{q-1}r_{vx}(A)) dv $$ Since $r_{x}$ is a morphism, $|\tau|$
is central, and $\bTr 1$ is a trace
\cite{Melrose46}, the tracial property of $\bTr q$ follows.   
  
To complete the proof, we need only prove that $\bTr q$ extends
$\aTr$, and this follows by integration in polar coordinates.
\end{proof}

It is not essential in the above statements that $A$ have integral
order.  Both the formula for $r_{x}(A)$ and the definition of $\bTr
q(A)$ make sense for $A \in \PsS {s}Y$, with $s$ not necessarily
integral.  We shall use this for operators of the form $D^{-z}A$ in
the next proposition. Actually, more is true of the $\bTr q$-traces of
elements of non-integral order than for elements of integral order:
let $s \not \in \ZZ$, then the action of $GL_q(\RR)$ by automorphisms
on $\PsS {s}{Z \times \RR^q}$ has the property
\begin{equation}\label{eq.act.GL}  
        \bTr q(T(A)) = |\det(T)|^{-1} \bTr q(A).  
\end{equation}  
 
This follows by using one of the several equivalent definitions of
$\bTr q(A)$, for $A$ of non-integral order, provided in the following
Lemma.

\begin{lemma}\label{lemma.coinc}\  
Let $A \in \PsS {s}{Z \times \RR^q}$, $s \not \in \ZZ$.  Then the
functions $Tr( \hat{D}(\tau)^{-z} \hat{A}(\tau) )$, and $\int_{\RR^q}
Tr( \hat{D}(\tau)^{-z} \hat{A}(\tau) ) d\tau$ are holomorphic for
$\Re(z - s) > d + q$ and extend to holomorphic functions on $\CC
\smallsetminus (\ZZ + s)$. At $z_0\not \in s +\ZZ$, these holomorphic
extensions satisfy
\begin{multline*} 
	\bTr q(D^{-z_0}A) = \left( \int_{\RR^q} Tr( \hat{D}(\tau)^{-z}
	\hat{A}(\tau)) d\tau \right) \vert_{z = z_00}\\ =
	\int_{S^{q-1}}\left ( \ITr_0^\infty \big[ x^{q-1}Tr \big(
	\hat{D}(x\tau)^{-z} \hat{A}(x\tau) \big) \big] \vert_{z = z_0}
	dx \right ) d\tau.
\end{multline*} 
\end{lemma}

\begin{proof}\ 
The function $\bTr q(D^{-z} A)$ is defined for all $z$, and it is seen
to be holomorphic on $\CC \smallsetminus (\ZZ + s)$ using the
definition of $\bTr q$ and Corollary \ref{cor.ext.one}, which also
gives the existence of the desired holomorphic extensions.

Finally, using again Corollary \ref{cor.ext.one} and Lemma
\ref{Lemma.on.trace}, we see that $\bTr q(D^{-z}A)$ and all the other
functions in the stated equation coincide for $\Re(z - s) > d +
q$. Because they are holomorphic on a connected open set containing
both $z_0 \in \CC \smallsetminus (s + \ZZ)$ and $\Re(z - s) > d + q$
($d = \dim Z$), they must coincide at $z_0$ also.
\end{proof}

\begin{proposition}\label{Prop.4}\ 
For any invertible, positive element $D\in\PsS {1}{Z \times \RR^q}$
and any holomorphic function $A: \CC \to \PsS m{Z \times \RR^q}$, the
function
\begin{equation*}  
        F_D(A;z)=\bTr q(D^{-z}A(z)), \quad \Re(z) > m + q + \dim Z,
\end{equation*}  
is holomorphic on $(\CC \smallsetminus \ZZ) \cup \{\Re z > m+q+\dim
Z\}$ with at most simple poles at the integers. The residue of this
holomorphic function at $0$ depends only on $A(0)$ and will be denoted
by $\RTr(A(0))$.  Moreover, $\RTr(A(0))$ vanishes on regularizing
elements, is independent of $D$, and defines a trace on $\PsS \infty{Z
\times \RR^q}$.
\end{proposition}

\begin{proof}\   
The function $F_D$ is holomorphic on the indicated domain by 
\ref{lemma.coinc}. Its poles are simple by Corollary \ref{cor.ext.one}
and by the definition of $\bTr q$ in terms of $\bTr 1$, see Equation
\eqref{eq.Tr.q}.
 
For the rest of the proof, it is enough to assume that $A(z)$ is
independent of $z$. We set then $A = A(z) = A(0)$. For example, the
proof of the fact that $\RTr$ is a trace and that it is independent of
the choice of $D$ is obtained from a standard reasoning, as
follows. We first write
\begin{equation*}   
        \bTr q(D^{-z}[A,B])=\bTr q(D^{-z}[D^{-z},D^{z}A]B)
\end{equation*}   
and observe that $[D^{-z},D^{z}A]B$ is a holomorphic function
vanishing at $0$.  This shows that $\bTr q$ vanishes on
commutators. The independence of $\bTr q$ on $D$ is a consequence of
\begin{equation*} 
        \bTr q(D^{-z}A)-\bTr q(D_1^{-z}A)=\bTr q(D^{-z}(\Id -D^{z}D_1^{-z})A),  
\end{equation*}  
using that $(\Id -D^{z}D_1^{-z}) A$ is a holomorphic function
vanishing at $0$.
\end{proof}

\begin{proposition}\label{Prop.reg}\ 
Let $D \in \PsS {1}{Z \times \RR^q}$, $D \ge 0$, invertible as above,
and let $A \in \PsS m{Z \times \RR^q}$. We denote $B_z = D^{-z}A \in
\PsS {m - z}{Z \times \RR^q}$, and write the asymptotic expansion  
\begin{equation}\label{eq.asymp.exp} 
        Tr(\hat{B}_z(\tau)) \sim \sum_{k \le m + d}  
        \beta_{k}(z,\tau/|\tau|)|\tau|^{k-z},\quad |\tau| \to \infty,  
\end{equation}  
defined for $\Re(z)$ large or for $z$ not an integer. Then the 
coefficient $\beta_{-1}$ is holomorphic in a neighborhood of $0$, 
        $\RTr (A) = \int_{| \tau | = 1 } \beta_{-1}(0,\tau),$  
and   
\begin{equation*} 
        \lim_{z \to 0} (\bTr q(D^{-z}A) - z^{-1}\RTr(A))  
        = \bTr q(A) + \int_{|\tau | = 1 }   
        \pa_z \beta_{-1}(z,\tau)\vert_{z = 0}.  
\end{equation*}  
\end{proposition}

\begin{proof}\ By the definition of $\bTr q$ using integration with  
respect to the orthogonal group, it is enough to prove the result for 
$q = 1$. By definition, $\bTr 1(A)$ is completely determined by 
$Tr(\pa_{x}^k \hat{A}(x))$.  This reduces our analysis to a lemma 
about integrals of functions in ${\mathcal M}_s$. 
  
Fix $\epsilon > 0$, and define ${\mathcal C}_\epsilon$ to be the space of  
functions $f(s,x)$, $- \epsilon < \Re (s) < \epsilon$, $ x \ge 0$, with  
the following properties:  
\begin{enumerate}  
\item {\em $f(s,x)$ is smooth in $(s,x)$ and holomorphic in $s$, for each   
fixed $x$};  
\item {\em For any $M \in \NN$ there exist $R >0$, $c_k \in \CC$, and 
complex valued functions $\alpha_k(s), b_k(s)$, and $h_M(s,x)$ 
satisfying 
$$  
        f(0,x)=h_M(0,x) + \sum_{k=-M}^{-1}\alpha_k(0) x^{k} + \sum_{k 
        = 0}^N \big(\alpha_k(0) + b_k(0) + c_k \log x \big)x^k, 
$$  
for $x \ge R$, and 
$$  
        f(s,x)=h_M(s,x)x^{-s} + \sum_{k=-M}^{N}\alpha_k(s) x^{k-s} + 
        \sum_{k = 0}^{N} \big( b_k(s) + c_k s^{-1} (1-x^{-s}) 
        \big)x^k, 
$$  
for $s \not = 0$, $x \ge R$, and $-\epsilon < \Re(s) < \epsilon$.} 
  
\item {\em $a_k(s)$, $b_k(s)$ are holomorphic in $s$, $h_M(s,x)$ is  
holomorphic in $s$, for each fixed $x$, and $h_M(s,\cdot) \in  
{\mathcal S}^{-M-1}([0,\infty))$, for each fixed $s$ in the strip  
$-\epsilon < \Re(s) < \epsilon$.}  
\end{enumerate}  
  
Of course, the choice of $R > 0$ is not important in the above definition.  
 
A simple but crucial observation is that $I(f) \in {\mathcal
C}_\epsilon$ if $f \in {\mathcal C}_\epsilon$.  Moreover, for any $f
\in {\mathcal C}_\epsilon$, with the coefficients of its canonical
asymptotic expansion denoted $\alpha_k$, $b_k$, and $c_k$, as in the
above definition, we have $$ \lim_{s \to 0} \ITr f(s,x)dx -
\alpha_{-1}(s)s^{-1}R^{-s} = \ITr f(0,x)dx + \alpha_{-1}(0)\log R, $$
which gives
\begin{equation}\label{eq.final}  
	\lim_{s \to 0} \big ( \ITr_0^\infty f(s,x)dx -
	\alpha_{-1}(0)s^{-1} \big ) = \ITr_0^\infty f(0,x)dx +
	\alpha_{-1}'(0).
\end{equation}  

The main idea is to prove that the familiar function 
\begin{equation*} 
	f_l (s,x) := Tr \big[( \pa_x^l (\hat{D}^{-s}(x) \hat{A}(x) +
	\hat{D}^{-s}(-x) \hat{A}(-x) ) \big ]
\end{equation*} 
is in $\mathcal C_\epsilon$, for $l$ large. Then $I^l(f_l) \in
\mathcal C_\epsilon$. This is, of course, a refinement of 
Proposition \ref{prop.ext}.
 
Now, by definition, $\bTr 1(D^{-s}A) =
\ITr_0^\infty I^l(f_l)(s,x) dx$, for any $l > m + d + 1$. For any 
such $l$, $f_l(s,x)$ is holomorphic in $s$, for $\Re(s) > -1$, and
hence $f_l \in \mathcal C_\epsilon$ (and $c_k = 0$, but this will play
no role in our reasoning), by Proposition
\ref{prop.ext}. Consequently, $I^l(f_l)(s,x) \in
\mathcal C_\epsilon$. Let $\alpha_k$ be the corresponding coefficients 
in the canonical asymptotic expansion of $I^l(f_l)$. 
 
We know, by classical results, that the functions $I^l(f_l)(s,x)$ and
$f_0(s,x)$ have holomorphic extensions in $s \not \in \ZZ$.  Moreover,
for $s \not \in \ZZ$, the difference $I^l(f_l)(s,x) - f_0(s,x)$ is a
polynomial in $x$, which shows that the coefficients of $x^{-1}$ in
the asymptotic expansions of these two functions are the
same. Consequently, $\beta_{-1}(s) = \alpha_{-1}(s)$, for $s \not =
0$. But by the definition of $\mathcal C_{\epsilon}$, $\alpha_{-1}$ is
holomorphic in a neighborhood of $0$, and hence $\beta_{-1}$ has a
holomorphic extension to a neighborhood of $0$, as claimed.
 
Finally, using $\bTr 1 (D^{-s}A) = \ITr_0^\infty I^l(f_l)(s,x)dx$, for $z$ in 
a small neighborhood of $0$,  and  
Equation \eqref{eq.final}, for $f=I^l(f_l)$, we obtain  
\begin{multline*} 
	\lim_{s \to 0} \big( \bTr 1(D^{-s}A) - \beta_{-1}(0)s^{-1}
	\big) = \lim_{s \to 0} \big( \ITr_0^\infty I^l(f_l)(s,x)dx -
	\alpha_{-1}(0)s^{-1} \big) \\ = \ITr_0^\infty I^l(f_l)(0,x) dx +
	\alpha_{-1}'(0) = \bTr 1(A) + \beta_{-1}'(0).
\end{multline*} 
This completes the proof. 
\end{proof}

We now drop the assumption that $B$ be reduced to a point.  To extend  
the above results to the general case, we proceed to a large extent as  
we did when $B$ was reduced to a point.  
   
Fix an invertible positive operator $D \in \PsS {1}{Y}$, and let $\CC  
\ni z \to A(z) \in \PsS {m}{Y}$ be an entire function.  Then  
\begin{equation}\label{eq.f.z.fam}  
        f_z(\tau)=Tr\big(\hat{D}(\tau)^{-z}|\tau|^k\hat{A}  
        (\tau)\big),\quad k \in \ZZ_+,  
\end{equation}   
is defined and holomorphic for any $\Re(z)> m+d = m + \dim Z$ and  
any fixed $\tau \in \GR^*$; moreover, the function $z \to f_z(\tau)$  
has a meromorphic extension to $\CC$, for each fixed $\tau$, with at  
most simple poles at integers.  
 
Let $\Omega = (\CC \smallsetminus \ZZ) \cup \{z,\Re(z)> m + \dim Z\}$.

\begin{lemma}\label{lemma.ext.fam}\  
Let $A(z) \in \PsS {m}{Y}$ be an entire function.  Also, let  
$f_z(\tau)$ be as above, with $z \in \Omega$ and $\tau \in \GR^*$.  
  
(i) The function $f_z(\tau)$ is in $\CI(\Omega \times \GR^*)$ and the   
map $z \to f_z(\tau)$ is holomorphic on $\Omega$, for each fixed $\tau  
\in \GR^*$.  
  
(ii) There is $g \in \CI(\Omega \times \GR^*)$ such that $g(z, \cdot) 
\in S^{m+d}(\GR^*)$, for each fixed $z$, such that $g(z,\tau)$ is 
holomorphic in $z$, for each fixed $\tau$, and such that 
$f_z(\tau)=g(z,\tau)|\tau|^{k-z}$, for all $\tau \ge 1$. Consequently, 
$f_z(x\tau) \in \mathcal{M}_{m+d+k-z}$, for all $\tau \not = 0$ 
($\,d=\dim Y - \dim \GR\,$). 
  
(iii) The function $z \to \ITr f_z$ is holomorphic on   
$\CC \smallsetminus \ZZ$, with at most simple poles at integers.  
\end{lemma}

\begin{proof}\ Since all the statements of the above theorem are  
statements about the local behavior in $b \in B$ of certain functions,  
we may assume that $Y$ is a flat $\GR$-space. This means, we recall,  
that $Y = B \times Z \times \RR^q$ and $\GR = B \times \RR^q$. Then  
we just repeat the proof of Proposition \ref{prop.ext} including an extra  
parameter $b \in B$, with respect to which all functions involved are  
smooth.  
\end{proof}

We now extend the definition of the   
various traces and functionals we considered above when  
$B$ was reduced to a point. This is not completely canonical, because  
we need to fix a metric on $\GR^*$ in order to obtain a volume form  
on the fibers of $\GR^* \to B$. The choice of the metric defines the  
group $O(\GR)$ of fiberwise orthogonal isomorphisms of $\GR$.  We also  
fix a lifting $Y/\GR \to Y$, which gives an isomorphism   
$$  
        Y \simeq Y/\GR \times_B \GR.  
$$   
This isomorphism and the metric on $\GR$ give an action of $O(\GR)$ on  
$Y$, which normalizes the structural action of $\GR$ by translations  
on $Y$. Consequently, the group $O(\GR)$ acts by isomorphisms on the  
algebra $\Psm {\infty}Y$. By Lemma \ref{Lemma.dilations}, the group  
$O(\GR)$ also acts by isomorphisms on $\PsS {\infty}Y$.

\begin{lemma}\  
Fix a metric on $\GR$ and choose a lifting $Y/\GR \to Y$, which gives rise  
then to an isomorphism $Y \simeq Y/\GR \times_B \GR$ and an action by  
automorphisms of the group $O(\GR)$ on $\PsS {\infty}Y$, as  
above. Then there exists an $O(\GR)$-linear map  
$$   
        E_Y : \PsS{\infty}Y \to \PsS {\infty}Y^{O(\GR)}   
$$    
such that $E_Y(AB) = AE_Y (B)$ and $E_Y (BA) = E_Y (B)A$, for all $A  
\in \PsS {\infty}Y^{O(\GR)}$ and $B \in \PsS {\infty}Y$.  Moreover,   
$\PsS {\infty}Y^{O(\GR)} \simeq \PsS{\infty}{Y/\GR \times  
\RR}^{\ZZ/2\ZZ},$ and hence the isomorphism class of the algebra  
$\PsS {\infty}Y^{O(\GR)}$ depends only on $Y/\GR$.  
\end{lemma}

\begin{proof}\ If $Y$ is a flat $\GR$ space, then this result follows  
right away from Lemma \ref{Lemma.restr.iso}. In general, we can choose  
trivializations of $Y$ such that the transition functions preserve the  
metric on $\GR$, and hence the transition functions are in $O(\GR)$.  
Because the isomorphisms of Lemma \ref{Lemma.restr.iso} commute with  
the action of the orthogonal group, the result follows.  
\end{proof}

We now consider $\CI(B)$--linear traces on $\PsS {m}Y$, for a general  
$\GR$ space $Y$. That is, we consider $\CI(B)$--linear maps $T: \PsS  
{m}Y \to \CI(B)$ such that  
$$   
        T(fA) =f T(A),\quad \text{for } f \in \CI(B) \text{ and } A  
        \in \PsS {m}Y,  
$$   
and    
$$   
        T([A,B]) = 0,\quad \text{for } A, B \in \PsS {m}Y.  
$$   
   
If we fix a metric on $\GR$, then we obtain a $\CI(B)$--linear trace  
$\bTr Y$ that generalizes the $\bTr q$--traces as follows. Suppose $Y  
= B\times Z \times \RR^q$ and $A =(A_b) \in \PsS {m}{Y}$, then we  
set  
\begin{equation}  
        \bTr Y (A)(b) = \bTr q (A_b).  
\end{equation}  
Because trace $\bTr Y$, for $Y = B \times Z \times \RR^q$, is  
invariant with respect to the action of the orthogonal group, the  
choice of an isomorphism $Y \cong Y/\GR \times_B \GR$ of $\GR$-spaces  
and of a metric on $\GR$ allow us to extend the definition of $\bTr Y$  
to arbitrary $Y$. We stress that this trace depends on the choices we  
have made. This new trace satisfies  
\begin{equation}  
        \bTr Y (A) = \bTr {Y/\GR \times \RR} (E_Y (A)).  
\end{equation}  
  
We then have the following immediate generalization of Proposition  
\ref{Prop.4} above:

\begin{proposition}\label{Prop.4.fam}\   
Fix a lifting $Y/\GR \to Y$, which is uded to define $\bTr Y (\,\cdot\,)$,  
as above. For any self-adjoint, invertible, positive element $D\in\PsS  
{1}{Y}$ and any holomorphic function $A: \CC \to \PsS m{Y}$, the  
function  
\begin{equation*}  
        F_D(A;z)=\bTr Y(D^{-z}A(z))  
\end{equation*}  
is holomorphic in $(\CC \smallsetminus \ZZ) \cup \{ \Re z> m+q+\dim Y 
-\dim \GR \}$ and has at most simple poles at the integers. The 
residue of this holomorphic function depends only on $A(0)$ and will 
be denoted by $\RT_Y(A(0))$.  Moreover, $\RT_Y(A(0))$ vanishes on 
regularizing elements, is independent of $D$, and defines a 
$\CI(B)$-linear trace on $\PsS m{Z \times \RR^q}$. This trace is 
independent of the choice of the isomorphism $Y \simeq Y/\GR \times_B 
\GR$.  
\end{proposition}

\begin{proof}\ Everything in this proposition follows from the  
case when $B$ is reduced to a point, except the independence of  
isomorphism $Y \simeq Y/\GR \times_B \GR$. For this we also use  
Equation \eqref{eq.act.GL}.  
\end{proof}

We also note that Proposition \ref{Prop.reg} extends virtually without  
change to families (that is, to the case $B$-nontrivial).  
   
The traces $\bTr q$ and $\bTr Y$ extend to matrix algebras by taking  
the sum of the traces of the entries on the main diagonal.

\section{Local index formulae\label{S.HI}}

We now return to the the study of the index of a family of invariant, 
elliptic operators. More precisely, we want local formulae for 
$Ch(\ind_a(A))$ when the fibers $\GR_b$ of $\GR \to B$ are 
simply-connected abelian Lie groups, that is, when $\GR$ is a {\em 
vector bundle}.  To this end, we shall use regularized traces and 
their properties developed in the previous section. 
  
If $A$ is a family of Dirac operators and $\GR$ is trivial, local 
formulae for $Ch(\ind_a(A))$ were obtained using heat kernels by 
Bismut in a remarkable paper, \cite{Bismut-fam}.  Our results are a 
step towards a similar result for arbitrary families of 
pseudodifferential operators invariant with respect to a bundle of Lie 
groups. The choice of the case when $\GR$ is a vector bundle may seem, 
at first sight, to be a long way from the general case. It is clear 
however that each type of Lie group bundle will have its specific 
features and hence have to be treated separately.  From the point of 
view of $K$-theory, vector bundles and bundles of simply-connected 
solvable Lie groups behave quite similarly. Moreover, by homotopy, the 
general case of simply-connected solvable Lie groups can be reduced to 
the particular case of a vector bundle, as we have shown in Section 
\ref{S.Homotopy.Inv}. 
  
No immediate generalizations of the results of this section to other 
classes of Lie group bundles seem possible. The case when the fibers 
of $\GR \to B$ are not simply-connected solvable has a quite different 
flavor, and even the case where $B$ is reduced to a point and $\GR$ 
has compact fibers is not well enough understood from the local 
perspective adopted in this section.  The case when $\GR \to B$ has 
simply-connected, non-commutative solvable fibers might seem more 
manageable at first sight, however, then we have some very difficult 
issues related to the choice of algebras closed under holomorphic 
functional calculus, not to mention that these groups may be not type 
I, so that the Fourier transform approach has no meaning. 
  
This should fully justify the choice of treating the case of vector 
bundles in such a great detail. 
  
Fix now a group morphism $\chi : H^{m}_c(\lgg^*) \to \CC$.  We want to 
obtain local formulae for $\chi(Ch(\ind_a(A)))$. For simplicity, we 
shall that $B$ and $Y/\GR$ are compact. Since the general case 
requires some additional ideas, we shall assume first that $Y$ is a 
flat $\GR$-space, that is, that $\GR = B \times \RR^q$ and $Y = B 
\times Z \times \RR^q$. 
  
Our approach is based on an interpretation of $\chi(Ch(\ind_a(A)))$  
using the Fedosov (or $\star$) product. We begin by recalling the  
definition of the Fedosov product and by making some general remarks  
on traces and their pairing with the $K$-theory of the algebras we  
consider.  
  
Let ${\mathfrak A}=\oplus_{k=0}^N {\mathfrak A}_k$, $N < \infty$ be a graded algebra  
endowed with a graded derivation $d:{\mathfrak A}_k \to {\mathfrak  
  A}_{k+1}$, the {\em Fedosov product} is  
defined by  
$$  
        a\star b = ab + (-1)^{\deg a}(da)(db).  
$$  
(The name is due to Cuntz and Quillen who have thoroughly studied the  
Fedosov product in connection to their approach to Non-commutative de  
Rham cohomology, see \cite{CuntzQuillen97}.)  We shall denote by  
$Q{\mathfrak A}$ the algebra ${\mathfrak A}$ with the Fedosov (or  
$\star$) product and by $Q_{ev}{\mathfrak A} \subset Q{\mathfrak A}$  
the subalgebra of even elements.  
  
Since we shall work with non-unital algebras also, it is sometimes  
necessary to adjoin a unit ``1'' to $Q{\mathfrak A}$.  The resulting  
algebra will be simply denoted by $Q^+{\mathfrak A} \simeq Q{\mathfrak  
A}\oplus \CC$. Similarly, $Q_{ev}^+{\mathfrak A} :=Q_{ev}{\mathfrak A}  
\oplus \CC$.  
  
A graded trace $\tau$ on $Q^+{\mathfrak A}$ restricts to an ordinary   
trace on $Q_{ev}^+{\mfk A}$, and hence it gives rise to a morphism  
$$  
        \tau_*:K^{\alg}_0(Q_{ev}^+{\mathfrak A}) \longrightarrow \CC,  
        \quad \tau_*[e] = \sum_j \tau(e_{jj}),  
$$  
for any idempotent $e = [e_{ij}] \in M_k(Q^+_{ev}\mfk A)$.  
If $\pi_*:K^{\alg}_0(Q{\mathfrak A})\to K^{\alg}_0({\mathfrak A}_0)$  
is the natural morphism induced by $\pi:Q_{ev}^+{\mathfrak A} \to  
{\mathfrak A}_0 \oplus \CC$, then $\pi_*$ is an isomorphism, by  
standard algebra results. Consequently, the trace $\tau$ also gives  
rise to a morphism  
\begin{equation}\label{eq.pairing.tilde}  
        \widetilde{\tau}:= \tau_* \circ \pi_*^{-1}:   
        K^{\alg}_0({\mathfrak A}_0 \oplus \CC)    
        \longrightarrow \CC,  
\end{equation}  
see \cite{ConnesNCG,Quillen}.  
  
The explicit form of the morphism $\widetilde{\tau}$ is not difficult  
to determine. Let $e\in {\mathfrak A}_0 \oplus \CC$ be an idempotent,  
then  
\begin{equation}  
        \bar{e}=\frac{1}{2} + \sum_{k\ge 0}^\infty
        (-1)^k\frac{(2k)!}{(k!)^2}\big( e - \frac{1}{2} + dede \big)
        (de)^{2k}
\end{equation}   
is an idempotent in $Q_{ev}^+{\mathfrak A}$ lifting $e$. (Note that 
the sum defining $\bar{e}$ is actually finite.) Assume the trace 
$\tau$ is concentrated on ${\mathfrak A}_{2k}$, $k \in 
\NN$ and $\tau(d{\mathfrak A}_{2k-1}) = 0$. Then the explicit formula  
for $\widetilde\tau([e])$ is 
\begin{equation}\label{eq.explicit.p}  
        \widetilde\tau([e]) = (-1)^k \frac{(2k)!}{k!} \tau((edede)^k). 
\end{equation}  
(We used $e(de)^k = (edede)^k$, valid for all $e$ satisfying $e^2 = e$.)  
  
Traces on $Q^+{\mfk A}$ are easy to obtain. Indeed, if $\tau$ is an  
even graded trace on ${\mathfrak A}$ satisfying $\tau({\mathfrak  
  A}_k)=0$, if $k \not = p$, and $\tau(d{\mathfrak A})=0$, then  
$$   
        \tau(a\star b - (-1)^{ij}b \star a) = 0,   
$$   
for any $a \in \mathfrak {\mathfrak A}_i$ and $b \in \mathfrak  
{\mathfrak A}_j$, and hence $\tau$ defines a graded trace on  
$Q{\mathfrak A}$. The trace $\tau$ defined above then extends to a  
trace on $Q_{ev}^+{\mathfrak A}$ by setting $\tau(1)=0$.  
  
We now define the algebras to which we shall apply the above  
considerations, if $Y = B \times Z \times \RR^q$ is a flat $\GR$  
space ($\GR = B \times \RR^q$). Let $\Omega^*(B)$ be the space of smooth forms on $B$,  
and consider the differential graded algebra  
\begin{multline*}  
        {\mathfrak A} := \PsS {\infty}{Y;\CC^N} \otimes_{\CI(B)}  
        \Omega^*(B) \otimes \Lambda^*\RR^q \simeq M_N(\PsS  
        {\infty}{Y}) \otimes_{\CI(B)} \Omega^*(B) \otimes  
        \Lambda^*\RR^q\\  
\simeq M_N(\Psi^{\infty}(Y)) \otimes \Omega^*(B) \otimes  
        \Lambda^*\RR^q \quad \text{ for } Y \text{ a flat } \GR-\text{space}.  
\end{multline*}  
Let $d_B$ the the de Rham differential in the $B$ variables (here we  
use the assumption that $Y$ is a flat $\GR$-space).  The differential  
of ${\mfk A}$ is then the usual de Rham differential $d_{DR}$:  
\begin{equation*}  
        d_{DR}(A)=\sum [t_i,A]d\tau_i + d_B(A),\; \; \Mand \; d(A  
        \xi)=d(A)\xi,  
\end{equation*}   
if $A \in \PsS {\infty}{Y;\CC^N} \otimes_{\CI(B)} \Omega^*(B) \simeq  
M_N(\PsS {\infty}Y) \otimes_{\CI(B)} \Omega^*(B)$ and $\xi$ is a  
product of some of the ``constant'' forms $d\tau_1, \ldots, d\tau_q$.  
Thus the differential $d_{DR}$ is with respect to the $B \times \RR^q$  
variables.  
  
Let $\pi^* \Lambda^* T^*B$ be the pull back to $Y$ of the exterior  
algebra of the cotangent bundle of $B$ and $F = \pi^* \Lambda^* T^*B  
\otimes \Lambda^*\RR^q$. Then  
$$  
        {\mathfrak A} \simeq \PsS {\infty}{Y;F \otimes \CC^N}.  
$$  
Inside ${\mathfrak A}$ we have the ideal of regularizing  
operators  
$$  
        {\mathfrak I} := \PsS {-\infty}{Y;\CC^N}  \otimes_{\CI(B)} \Omega^*(B)  
        \otimes \Lambda^*\RR^q \simeq \PsS {-\infty}{Y;F \otimes \CC^N},  
$$   
with quotient algebra  
\begin{multline*}  
        \mathfrak B := \mathfrak A/\mathfrak I = M_N(\PsS  
        {\infty}{Y}/\PsS {-\infty}{Y}) \otimes_{\CI(B)}  
        \Omega^*(B)\otimes \Lambda^*\RR^q \\ \simeq \PsS {\infty}{Y;F  
        \otimes \CC^N} / \PsS {-\infty}{Y;F \otimes \CC^N}.  
\end{multline*}  
   
We consequently obtain the exact sequence of algebras  
\begin{equation}\label{eq.exact.F}  
        0 \to Q_{ev}{\mathfrak I} \to Q_{ev}{\mathfrak A} \to  
        Q_{ev}{\mathfrak B} \to 0,  
\end{equation}   
which gives rise to the boundary map $\pa_Q$,  
$$   
        \pa_Q : K_1^{alg}(Q_{ev}{\mathfrak B}) \to  
        K_0^{alg}(Q_{ev}{\mathfrak I}) =  
        K_0^{alg}(Q_{ev}^+{\mathfrak I}_0),  
$$   
on algebraic $K$-theory.  
  
We now define the traces we shall consider on the algebras $\mfk 
A,\mfk B, \ldots, Q\mfk I$. Let $\omega : \Omega^{k}(B) \to \CC$ be a 
closed current, that is, a continuous map such that $\omega(d\eta) = 
0$ for any form $\eta$.  Then $\omega$ defines a morphism $\chi_\omega 
:H^*_c (\lgg^*) \to \CC$.  Similarly, let 
\begin{equation*}  
        A \otimes \zeta \otimes \xi \in {\mathfrak A} := \PsS 
        {\infty}{Y;\CC^N} \otimes_{\CI(B)} \Omega^*(B) \otimes 
        \Lambda^*\RR^q, 
\end{equation*}  
and define  
\begin{equation}\label{eq.def.tau.o}  
        \tau_\omega(A \otimes \zeta \otimes \xi) = \bTr q (A) 
        \omega(\zeta) \int_{\RR^q} \xi. 
\end{equation}  
Then $\tau_\omega$ is a trace on ${\mathfrak A}$ satisfying  
$\tau_{\omega}(d(a))= 0$ if $a \in {\mathfrak I}$.    
  
If $e = (e_0, \lambda) \in M_N(\PsS {-\infty}Y + \CC) = M_N(\PsS 
{-\infty}Y) \oplus M_N(\CC)$ is an idempotent, then $e$ defines a 
class $x_e \in K^0(\lgg^*) \cong K^0(B \times \RR^q)$.  It is useful 
to review the definition of this class.  First, we can replace 
$e=(e_0,\lambda)$ by an equivalent projection such that there exists 
another projection $p=(p_0,\lambda)$ satisfying $dp = 0$ and $pe = ep 
= p$. More precisely, the projection $p$ is such that its indicial 
family, $\hat{p}:=\hat{p_0} + \lambda$ consists of a smooth family of 
projections on $\GR^*$, acting on the fibers of $\GR^*\times_B Y \to 
\GR^*$, with values in $M_N(\Psi^{-\infty}(Y_b) + \CC)$, and constant  
along the fibers of $\GR^* \to B$. This projection $p$ has the  
property that $p {\mfk I}_0 \subset {\mfk I}_0$. Then we define the  
vector bundle $V_e$ on $\lgg^*$ such that its fiber at $\tau$ is the  
range of the indicial operator $\hat{e}(\tau) - \hat{p}$.  Finally,  
the class $x_e$ defined by $e$, which we are looking for, is  
$[V_e]-r[1]$, where $r$ is the rank of $V_e$.  
  
It is interesting to compare the Chern character of the bundle $V_e$  
to the pairing $\widetilde{\tau}_\omega[e]$. First, the vector bundle $V_e$  
is trivial at infinity. The curvature of the Grasmannian connection  
$\nabla^e = e\circ d$ is  $R^e := (e \circ d)^2 = edede$, by a  
standard computation.  If $k +q = 2p > 0$ is even, then  
$\tilde\tau_\omega$ is an even, graded trace, and hence it pairs with  
$[e]$. Using the explicit formula from Equation \eqref{eq.explicit.p},  
we obtain that  
\begin{equation}\label{eq.30}   
        \tilde{\tau}_{\omega}[e] = \frac{(2p)!}{p!}\chi_\omega\big 
        ((-R^e)^p/p! \big) 
\end{equation}   
recovers (up to a constant) the pairing of $Ch(V_e)$, the Chern 
character of $V_e$, with the cohomology class of $\chi_\omega$. In 
the notation introduced above, we have 
\begin{equation}\label{eq.pair.K}   
        \tilde\tau_{\omega}[e] = 
        \frac{(2p)!}{p!}\chi_{\omega}(Ch(V_e)) = 
        \frac{(2p)!}{p!}\chi_\omega(Ch(x_e)). 
\end{equation}   
See \cite{ConnesNCG}. It also follows that all morphisms  
$K_0^{alg}(\PsS {-\infty}Y) \simeq K_0(\PsS {-\infty}Y) \to \CC$ are  
of the form $\tilde{\tau}_\omega$, for a suitable $\omega$. (Recall   
that $\tilde\tau_{\omega}$ is defined by Equations \eqref{eq.pairing.tilde}   
and \eqref{eq.explicit.p}.)  
  
Recall now that we defined the analytic index $\ind_a$ to be the  
composite map  
$$   
        \ind_a : K_1^{alg}({\mathfrak B}_0) \stackrel{\pa}{\longrightarrow}   
         K_0^{alg}({\mathfrak I}_0) \simeq  K^0(\lgg^*),   
$$   
where ${\mathfrak I}_0 = M_N(\PsS {-\infty}Y)$ and $\pa$ is the boundary  
map in algebraic $K$-theory associated to the exact sequence  
$0 \to {\mfk I}_0 \to {\mfk A}_0 \to {\mfk B}_0 \to 0.$ Moreover, for the case  
we are currently discussing, that where   
$Y$ is a flat $\GR$-space, $\lgg^* \cong B \times \RR^q$.  
Let $\chi :H^*_c (\lgg^*) \to \CC$ be an arbitrary group morphism.  
As we explained at the beginning of the section, we are interested in  
understanding the morphism   
$$  
        \chi \circ Ch \circ \ind_a : K^{alg}_1(\PsS {\infty}{Y;\CC^N}  
        / \PsS {-\infty}{Y; \CC^N}) =K^{alg}_1({\mfk B}_0) \to  
        \CC.  
$$  
By the above discussion and linearity, we may assume that  
$\chi=\chi_\omega$, for some closed current $\omega:\Omega^k(B) \to  
\CC$. Then a preliminary formula for the composition $\chi_\omega  
\circ Ch \circ \ind_a$ is given in the following lemma.

\begin{lemma}\label{Lemma.morphism}\  
Let $\omega: \Omega^k(B) \to \CC$ be a closed current such that $k + q  
=2p >0$ is even. Denote by $\chi_\omega : H^*_c(\lgg^*) \to \CC$ the  
morphism defined by $\omega$. Then  
$$  
        \chi_\omega \circ Ch \circ \ind_a = \tilde\tau_\omega  
        \circ \pa : K_1^{alg}({\mfk B}_0) \to \CC.  
$$  
\end{lemma}

\begin{proof}\ This follows by applying the above constructions to   
$\ind_a(A)$, $A$ elliptic, using Equation \eqref{eq.pair.K}.  
\end{proof}

We now turn to the computation of $\tilde\tau_{\omega}\circ \pa :   
K_1^{alg}({\mfk B}_0) \to \CC$.  
We shall use the generic notation $\pi$ for all the quotient morphisms  
$Q{\mathfrak A} \to {\mathfrak A}_0$ and $Q{\mathfrak B} \to  
{\mathfrak B}_0$, and $Q{\mathfrak I} \to {\mathfrak I}_0$.   
Also, we shall denote by  
$$   
        [A,B]_\star = A\star B - (-1)^{ij}B \star A,  
$$   
for $A \in {\mfk A}_i$, $B \in {\mfk A}_j$, the graded commutator in  
$Q{\mfk A}$ with respect to the $\star$-product. (Although in most  
cases $i$ and $j$ are even, as it is the case in the lemma below.)

\begin{lemma}\label{lemma.ind}\   
Let $u \in {\mfk B}_0$ be an invertible element with inverse  
$v$. Choose liftings $A$ and $B$ of $u$ and, respectively, $v$. Also,  
let $B' :=\sum_{k=0}^{\infty} (-1)^kB(dAdB)^k$ and let $\tau$ be a  
closed graded trace on ${\mfk I}$ satisfying $\tau([\mfk A, \mfk I]) =  
0$. Then  
\begin{equation*}   
        \widetilde\tau\circ \pa[u] = \tau\big ([A,B']_\star \big).  
\end{equation*}   
\end{lemma}

\begin{proof}\   
The map   
$$  
        \pi_*:K^{\alg}_1(Q{\mathfrak B}) \to K^{\alg}_1({\mathfrak B}_0)  
$$   
is onto because if $u \in {\mathfrak B}_0 = \PsS{\infty}{Y;\CC^N} /  
\PsS{-\infty}{Y;\CC^N}$ is invertible in ${\mathfrak B}_0$ with  
inverse $v$, then its image $u'$ in $Q{\mathfrak B}$ is also  
invertible with inverse  
$$  
        v':=\sum_{k=0}^{\infty} (-1)^kv(dudv)^k.  
$$   
(The sum is actually finite for our algebras.) The relations    
$$  
        u'\star v' = 1 = v' \star u'  
$$   
are easily checked. From the naturality of the boundary map in  
algebraic $K$--theory, we obtain that  
$$   
        \pi_*\circ \pa_Q =\pa \circ \pi_*:   
        K^{\alg}_1(Q {\mathfrak B}) \to K^{\alg}_0({\mathfrak I}_0),   
$$   
and hence   
\begin{equation}  
        \tilde\tau\circ \pa [u] = \tau_* \circ \pa_Q [u'].  
\end{equation}   
This simple relation will play an important role in what follows because   
it reduces the computation of $\tilde\tau \circ \pa$ to the   
computation of $\tau_* \circ \pa_Q$.   
   
Lift $u \in M_N(\mfk B_0)$ to an element $A \in M_N(\mfk A_0) =  
M_N(\PsS {\infty}{Y})$ and its inverse $v$ to an element $B \in M_N(\PsS  
{\infty}{Y})$, as in the statement of the lemma.  This gives for $v'$  
(the inverse of the image $u'$ of $u$ in $Q{\mfk B}$ with respect to  
the $\star$ product) the explicit lift   
$$  
        B' :=\sum_{k=0}^{\infty}(-1)^kB(dAdB)^k.  
$$  
We now proceed by direct computation (as in \cite{Nistor2}, for example), using the   
explicit formula $\pa_Q([u'])=[e_1]-[e_0]$ with $e_0=1 \oplus 0$ and  
\begin{equation}\label{eq.explicit.e1}  
	e_1=\left[\begin{array}{cc} 2A\star B'-(A\star B')^2 &   
	A\star (2-B'\star A)\star (1-B\star A)  
	\\(1-B'\star A)\star B' & (1-B'\star A)^2  
\end{array}\right].  
\end{equation}  
(all products and powers are with respect to the $\star$-product). Then,   
\begin{multline*}  
      \tau_*\circ \pa_Q([u'])= \tau([e_1]-[e_0]) = \\ \tau (2A\star B' 
      - (A\star B')^2 + (1 - B'\star A)^2 -1) = \tau( (1 - B'\star 
      A)^2 - (1 - A\star B')^2). 
\end{multline*}  
Then we notice that $\tau([A, B' - B'\star A \star B']_*) = 0$,  
because $B' - B'\star A \star B' \in \mfk I$. This relation and its  
analogue obtained by switching $A$ with $B'$ then give that  
$\tau_*\circ \pa_Q([u']) = \tau(A\star B'-B'\star A)$, and the lemma  
follows.  
\end{proof}

Let $D \in \PsS {1}{Y}$ be the operator used to define $\bTr q$ and  
$\RTr$ in the previous section. Also, let $\iota:\Lambda^q \RR^q \to  
\CC$ be the isomorphism given by contraction with the (dual) of the  
top form on $\RR^q$. This gives rise to maps  
\begin{equation}  
        \bTr Y \otimes \iota, \RT_Y \otimes \iota : Q\PsS {\infty}Y  
        \to \Omega^*(B),  
\end{equation}  
which vanish on forms of degree less than $q$ in $d\tau_1, \ldots,  
d\tau_q$.  Because $\bTr Y$ and $\RT_Y$ are $\CI(B)$-linear graded  
traces on ${\mfk A}_0$, $\bTr Y \otimes \iota, \RT_Y \otimes \iota$  
are $\Omega^*(B)$-linear (graded) traces.  If $\omega: \Omega^*(B) \to  
\CC$ is a closed current, we denote  
\begin{equation}\label{eq.def.ro}   
        \rho_\omega (A) = \langle \omega, \RT_Y \otimes \iota  
        (A) \rangle.  
\end{equation}   
and note that $\rho_\omega$ is a closed graded trace. Also, note that  
$\tau_\omega(A) = \langle \omega , \bTr Y \otimes \iota (A) \rangle$,  
which we shall use to extend $\tau_\omega$ to operators of  
non-integral orders, still preserving the tracial property. Moreover,  
$\tau_\omega(dA)=0$, if $A$ has non-integer order. Consequently,  
$\tau_\omega([A,B]_\star)=0$ if $\ord A + \ord B$ is not an integer.  
This is seen by noticing that $\tau_\omega( d(D^{-z}A ) = 0$ for  
$\Re(z)$ large first, and then for $z$ such that $z \ord D + \ord A$  
not an integer, by analytic continuation.   
  
We then have

\begin{lemma}\label{lemma.res}\ If $A(z) \in Q{\mfk A}$ is holomorphic in a   
neighborhood of $0 \in \CC \smallsetminus \ZZ^*$, then the function  
$\tau_\omega(D^{-z}\star A(z))$ holomorphic and at $0$ it has a simple pole   
with residue $\rho_\omega (A(0))$.  
\end{lemma}

\begin{proof}\ We have that   
$$  
        D^{-z} \star A(z) = D^{-z}A(z) + dD^{-z}dA(z).  
$$  
Now we observe that   
$$    
        \lim_{z \to 0} z^{-1} d D^{-z}= \lim_{z \to 0}     
        z^{-1} \sum  [t_j,D^{-z}]d \tau_j=   
        -\sum  [t_j,\log D]d \tau_j=-d\log D   
$$   
Consequently, $D^{-z}\star A(z) = D^{-z}B(z)$ for some holomorphic  
function $B$ such that $B(0) = A(0)$.  The result then is an immediate  
consequence of Proposition \ref{Prop.4.fam}.  
\end{proof}

Let $\pa : K_1^{alg}(\mfk B_0) = K_1^{\alg}(\PsS{\infty}Y / \PsS 
{-\infty}{Y})) \to K_0^{\alg}(\PsS{-\infty}Y)=K_0^{alg}(\mfk I_0)$ be 
the boundary map in algebraic $K$-theory, as above. Recall that the 
map $\tilde\tau_\omega : K_0^{\alg}(\mfk I_0) \to \CC$ is given by the 
Equations \eqref{eq.pairing.tilde} and \eqref{eq.30} and that 
$\rho_\omega$ is given by the Equation \eqref{eq.def.ro}.   
 
We continue to assume that $Y = B \times Z \times \RR^q$ is a flat 
$\GR = B \times \RR^q$-space, that $q > 0$, and that $B$ and $Z$ are 
compact.

\begin{theorem}\label{Theorem.Index.F}\   
Let $u \in M_N(\PsS {\infty}{Y}/ \PsS {-\infty}{Y})$ be an invertible 
element, and choose $A,B \in M_N(\PsS {\infty}{Y})$ such that $A$ maps 
to $u$ and $B$ maps to $u^{-1}$. If $\omega : \Omega^{k}(B) \to \CC$ 
is a closed current such the $k + q = 2p > 0$ is even, then 
\begin{multline*}  
  	\tilde\tau_\omega \circ \pa [u]= - 2(-1)^{p} \tau_\omega 
  	((dAdB)^{p}) = 2(-1)^{p} \tau_\omega ((dBdA)^{p}) \\ = 
  	-\rho_\omega \big ( u^{-1}[\log D, u](u^{-1}du)^{2p} \big ) 
  	-\rho_\omega \big ( u^{-1}[d\log D, du](u^{-1}du)^{2p-2} \big 
  	). 
\end{multline*}   
\end{theorem}

\begin{proof}\  We shall denote    
$$   
        {\mfk A} = M_N(\PsS {\infty}Y \otimes_{\CI(B)} \Omega^*(B) \otimes  
        \Lambda^*\RR^q),  
$$   
as before. Moreover, $\mfk B$ and $\mfk I$ will have the meaning they had  
before.  
  
We shall use Lemma \ref{lemma.ind}.  Let  
$B'= \sum_{k=0}^{\infty} (-1)^k B(dAdB)^k$ be as in that lemma and  
evaluate the commutator $[A,B']_\star$ (with respect to the $\star$  
product).  We obtain  
\begin{multline*}  
        [A,B']_\star = \sum_{l=0}^{\infty} (-1)^l AB(dAdB)^l   
        - \sum_{l=0}^{\infty} (-1)^l B (dAdB)^l A  \\  
        + \sum_{l=0}^{\infty} (-1)^l (dAdB)^{l+1}   
        - \sum_{l=0}^{\infty} (-1)^l (dBdA)^{l+1},  
\end{multline*}   
the sums being of course finite.  
  
We next observe that $\tau_\omega(AB(dAdB)^{l})  
=\tau_\omega(B(dAdB)^{l}A)$ and $\tau_\omega((dAdB)^{l})=  
-\tau_\omega((dBdA)^{l})$ because $\tau_\omega$ is a graded trace on  
$\mfk A$ (with the usual product).  Using this we obtain from Lemma  
\ref{lemma.ind} that  
\begin{equation}   
        \tilde\tau_\omega \circ \pa [u]=    
        \tau_\omega (A \star B' - B' \star A)=     
        2(-1)^{p-1}\tau_\omega ((dAdB)^{p}).  
\end{equation}  
This proves the first part of our formula.   
  
We now prove the second part of our formula.   
The commutator $[A,B']\star = A \star B' - B' \star A$ maps to $u\star  
v -v\star u = 0$ in $Q{\mfk B}$, and hence $[A,B']_\star$ is in  
$Q{\mfk I}$. Consequently,   
$$  
	\tau_\omega([A,B']_\star) = \lim_{z \to 0}  
	\tau_\omega(D^{-z}\star [A,B']_\star).  
$$  
Next we observe that $ \lim_{z \to 0} z^{-1} d D^{-z}= -d\log D $   
and hence $d\log D$ (defined in the canonical representation)  
is actually in $\PsS {-1}Y$, in spite of the fact   
that $\log D$ is not in the algebra $\PsS {\infty}Y$. Moreover, $z^{-1}dD^{-z}$   
is holomorphic at $0$.   
  
Using that $\tau_\omega([A,D^{-z}B']_\star) = 0$ for all $z$ such that 
$-z + \ord A + \ord B$ is not an integer, we finally obtain 
\begin{multline*}   
        \tau_\omega([A,B']_*) = \lim_{z \to 0}    
        \tau_\omega(D^{-z}\star [A,B']_\star) \\   
        = \lim_{z \to 0} \tau_\omega([D^{-z},A]_{\star} \star B')   
        = \lim_{z \to 0} z\tau_\omega(D^{-z} \star F(z)),   
\end{multline*}   
where $F(z) = z^{-1}[D^{-z},D^z \star A]_{\star} \star B'$. Since $F$ 
is a holomorphic function in a neighborhood of $0$ with  
\begin{equation*} 
	F(0) = -[\log D, A]_*B' - d[\log D,A]_*dB', 
\end{equation*} 
it further follows from Lemma \ref{lemma.res} that 
\begin{equation*}   
  	\lim_{z \to 0} z\tau_\omega(D^{-z} \star F(z)) = \rho_\omega  
  	(F(0)) = - \sum_k \rho_\omega \big ( u^{-1}[\log D,  
  	u]_*(u^{-1}du)^{2k}\big ),  
\end{equation*}   
because $\rho_\omega$ is a closed graded trace.  
  
Putting toghether the above formulae, we obtain the desired result.  
%
%
\end{proof}

The method we used in the previous theorem works even when $Y$ is not 
a flat $\GR$ bundle, with only minor changes. Actually, using a trick 
of Connes, we can formally treat the general case to a large extent as 
if it was a flat bundle.  For example, the definition of $\mfk A$, 
$\mfk I$, and $\mfk B$ extend to this case without change. However, 
$\omega$ and the differential $d$ cannot be defined as before to enjoy 
all the previous properties. Here is what changes. First, we need to 
take care of the fact that $\lgg$ may not be orientable. This is easy 
dealt with by considering linear functionals on forms twisted with the 
orientation sheaf instead currents (which are linear functional on 
ordinary forms).  We say then that the current $\omega$ must be 
twisted with the orientation sheaf of $\lgg$. Note that it still makes 
sense to talk about a closed twisted current, because the de Rham 
differential extends to the space $\Omega^j \otimes {\mathcal O}$ of 
forms twisted with the orientation sheaf:  
\begin{equation*} 
	d: \Omega^j \otimes {\mathcal O} \to \Omega^{j + 1} \otimes 
	{\mathcal O}, 
\end{equation*} 
and $d^2 = 0$. Once we understand the nature of $\omega$, the 
definitions of $\tau_\omega$ and $\rho_\omega$ carry through without 
any change. 
  
The problem is to define the differential structure on $\mfk A$ in the non-flat case.    
Choose a trivializing covering of $B$ and a partition of unity $\phi_\alpha^2$,   
subordinated to that covering. Using the construction from the flat case,  
on each of the trivializing open sets $U_\alpha$   
of the covering we have an operator $d_\alpha : \mfk A\vert_{U_\alpha}  
\to \mfk A\vert_{U_\alpha}$ (which satisfies $d_\alpha^2 = 0$, but this is   
irrelevant to us). Let $\nabla := \sum \phi_\alpha d_\alpha \phi_\alpha$.  
Then $\nabla : \mfk A \to \mfk A$ is a degree one derivation and   
$\nabla^2(A) = [\Theta, A]$, for some $\Theta \in {\mfk A}$. Moreover,  
$\tau_\omega(\nabla(a)) = 0$ if $a \in \mfk I$, $\rho_\omega(\nabla(a)) = 0$   
if $a \in \mfk I$, and both $\tau_\omega$ and $\rho_\omega$ are  
traces on $\mfk A$.  
  
The fact that $\nabla^2$ is not zero in general means that we cannot  
use it to define the Fedosov product using $\nabla$. This is not a big  
issue, however, because we can proceed as in \cite{ConnesNCG} (see  
also \cite{NistorSUPER}).  The idea, due to Connes, is to enlarge the  
algebra $\mfk A$ and to perturb $\nabla$ such that it becomes a  
differential. We now review this construction following  
\cite{NistorSUPER}.  
  
We first introduce a formal variable $X$  such that   
$$  
	a X b = 0 \quad \text{and} \;(aX)(Xb) = a\Theta b  
	\quad \text{if} \; a , b \in \mfk A  
$$   
(we never consider $X$ alone, only in formulae like $aX$, $Xb$, or $XaX$).   
Define $\overline{\mfk A} = \mfk A + X \mfk A + \mfk A X + X \mfk A X$, with  
the induced grading such that $\deg X = 1$. We define  
then $\overline{d}a = \nabla(a) + Xa + (-1)^{\deg a} aX$. We can then extend   
$\overline{d}$ to  
a derivation of $\overline{\mfk A}$ such that $\overline{d}(X) = 0$ and we can also  
extend $\tau_\omega$ to a graded trace, call it $\overline{\tau}_\omega$,  
on $\overline{\mfk A}$. The explicit formula for this trace is  
$$  
	\overline{\tau}_\omega(a_{00}+a_{01}X+Xa_{10}+Xa_{11}X) = \tau_\omega(a_{00}) -   
	(-1)^{\deg a_{11}} \tau_\omega (\Theta a_{11}).  
$$  
The preceeding theorem will then extend to this new  
setting where $\overline{\mfk A}$ replaces $\mfk A$, $\overline{d}$ replaces  
$d$, and $\overline{\tau}_\omega$ replaces $\tau_\omega$. The place of  
$\rho_\omega$ will be taken by $\overline{\rho}_\omega$, defined by a formula  
similar for that for $\overline{\tau}_\omega$ above:  
$$  
	\overline{\rho}_\omega(a_{00} + a_{01}X + Xa_{10} + Xa_{11}X) =   
	\rho_\omega(a_{00}) - (-1)^{\deg a_{11}} \rho_\omega (\Theta a_{11}).  
$$  
This shows that $\overline{\rho}_\omega$, unlike $\overline{\tau}_\omega$,   
is a {\em closed} graded trace ($\overline{\tau}_\omega$ is not closed).  
  
Note, in the following theorem, that the meaning of $\tilde \tau_\omega$  
does not change in the case of non-flat bundles.

\begin{theorem}\label{Theorem.Index.F'}\   
Let $u \in M_N(\PsS {\infty}{Y}/ \PsS {-\infty}{Y})$  and   
$A,B \in M_N(\PsS {\infty}{Y})$ be as in Theorem \ref{Theorem.Index.F'}.  
Also, let $\omega:\Omega^{k}(B) \to \CC$ is a closed twisted current (\ie with values   
in the orientation sheaf of $\lgg$) such the $k + q = 2p$ is even. If   
$\overline{d}$, $\overline{\rho}_\omega$, and $\overline{\tau}_\omega$ are as   
above, then  
\begin{multline*}  
  	\tilde\tau_\omega \circ \pa [u]= - 2(-1)^{p} \overline{\tau}_\omega   
	((dAdB)^{p} ) = 2(-1)^{p} \overline{\tau}_\omega ((dBdA)^{p})  \\  
  	= -\overline{\rho}_\omega \big ( u^{-1}[\log D, u](u^{-1}du)^{2p} \big )  
	-\overline{\rho}_\omega \big ( u^{-1}[d\log D, du](u^{-1}du)^{2p-2} \big ).  
\end{multline*}   
\end{theorem}

\begin{proof}\ The proof is word for word the same, if we replace  
$\tau_\omega$ with $\overline{\tau}_\omega$, $\rho_\omega$ with $\overline{\rho}_\omega$,   
$\mfk A$ with $\overline{\mfk A}$, and so on.  
\end{proof}

After the above discussion, the matter of extending 
Theorem~\ref{Theorem.Index.F} to the non-flat case seems a 
triviality. It is however far from being so, because the 
multiplication in the algebra $\overline{\mfk A}$ is more 
complicated. In particular, it introduces the curvature $\Theta$ of 
the vertical bundle $T_{vert}(Y/\GR)$. The appearance of $\Theta$ in 
formulae is actually a good thing, because we know from 
Theorem~\ref{Theorem.Chern} that the curvature has to be there. 
Nevertheless, this does not mean that we can easily deduce Theorem 
\ref{Theorem.Chern} from the Theorem~\ref{Theorem.Index.F'} 
above. Obtaining Theorem~\ref{Theorem.Chern} from the 
Theorem~\ref{Theorem.Index.F'} is a question that remains to be 
solved.

\section{Higher eta invariants in algebraic $K$-theory\label{S.HEF}}

We consider the same setting as in the previous section.   
More precisely,  $Z \to B$ is a fiber bundle,   
$\GR \to B$ is a vector bundle and $Y = Z \times_B \GR$, with the induced  
action of $\GR$. We state our results below for the case when  
both $Z$ and $\GR$ are trivial bundles. When dealing with non-trivial bundles,   
we replace $\tau_\omega$ with $\overline{\tau}_\omega$,  
$d$ with $\overline{d}$, and $\mfk A$ with $\overline{\mfk A}$,  
as we did at the end of that section.   
  
Consider then a closed current $\omega : \Omega^k (B) \to \CC$ and let  
$\tau_\omega$ be the associated trace on $\mfk A$. Then a direct  
computation gives that  
$$  
        \phi_\omega(a_0,a_1,\ldots ,a_l)=\tau_\omega(a_0da_1\ldots  
        da_l)/l!, \quad l = q + k  
$$  
is a $l$-Hochschild cocycle on $\PsS{\infty}Y$.   
The Dennis trace map \cite{Karoubi}  
$$  
        K_l^{\alg}(\PsS {\infty}Y) \to \Hd_l(\PsS {\infty}Y)  
$$  
and the morphism $\Hd_l(\PsS {\infty}Y) \to \CC$ defined by  
$\phi_\omega$ give rise by composition to a morphism  
\begin{equation}  
        \eta_\omega: K_l^{\alg}(\PsS {\infty}Y) \to \CC.  
\end{equation}  
  
Because the restriction of $\phi_\omega$ to $\PsS {-\infty}Y$ is  
cyclic (so it defines a cyclic cocycle), the composition  
$$  
        K_l^{\alg}(\PsS {-\infty}Y) \longrightarrow K_l^{\alg}(\PsS  
        {\infty}Y) \stackrel{\eta_\omega}{\longrightarrow} \CC  
$$  
factors as  
$$  
        K_l^{\alg}(\PsS {-\infty}Y) \longrightarrow K_l^{\topo}(\PsS  
        {-\infty}Y) \stackrel{(\phi_\omega)_*}{\longrightarrow} \CC,  
$$  
where $(\phi_\omega)_*$ is the pairing of cyclic homology with  
topological $K$-theory. In particular, $\eta_\omega$ is non-zero if  
$\omega$ is not exact.  
  
The morphism $\eta_\omega$ does not factor through topological $K$-theory  
though. This is seen by noticing that  
$$  
        K_l^{\topo}(\PsS {-\infty}Y)   
        \longrightarrow K_l^{\topo}(\PsS {0}Y)  
$$   
vanishes for any $p$, as proved in \cite{Melrose-Nistor1}. Moreover,  
for $B$ reduced to a point, $\GR = \RR$, $k =1$, $\omega (fdt) =  
\int_\RR f(t) dt$, and  
$$  
        D=D_0 + \pa_{t} \in M_N(\PsS {\infty}Y),  
$$   
the indicial map of an admissible (chiral) Dirac operator on $Y\times  
\RR$, the main result of \cite{Melrose46} states that  
$\eta_\omega(D)=\eta(D_0)/2$, where $\eta(D_0)$ is the  
``eta''--invariant introduced by Atiyah, Patodi, and Singer in  
\cite{APS}. As proved in \cite{Melrose46}, this gives that the  
eta invariant of $D_0$ is the value at $D$ of a group morphism  
$K_1(\PsS{\infty}{M \times \RR}) \to \CC$. This group morphism  
coincides with $\eta_\omega$, $\omega(f) = \int_\RR f$, in the above  
notation.

It is tempting then to try to define a higher eta invariant on 
$\PsS{\infty}{M \times \RR^q}$, $q = 2k -1$, by the formula 
$\eta_k(D_0) = \bTr q( (D^{-1} dD)^{2k-1})$, where $D = D_0 + 
c(\tau)$; however, as it was proved by Lesch and Pflaum, this is not 
multiplicative, and besides, it coincides with usual eta invariant of 
$D_0$ (up to a multiple depending only on $k$).

\providecommand{\bysame}{\leavevmode\hbox to3em{\hrulefill}\thinspace}

\end{document}